\documentclass[MSNbibl,number,citesort,dvips]{arxbj}
\usepackage{upgreek}

\aid{0}
\volume{20}
\issue{1}
\pubyear{2014}
\firstpage{1}
\lastpage{27}
\doi{10.3150/12-BEJ474} 

\makeatletter
\renewcommand{\pi}{\uppi}

\newcommand{\eqref}[1]{(\ref{#1})}
\newtheorem{theorem}{Theorem}
\newtheorem{proposition}[theorem]{Proposition}
\newtheorem{lemma}[theorem]{Lemma}
\newtheorem{corollary}[theorem]{Corollary}
\newproclaim{definition}[theorem]{Definition}
\newremark{remark}[theorem]{Remark}

\newremark{example}{Example}

\newcommand{\R}{\mathbb{R}}
\newcommand{\N}{\mathbb{N}}
\newcommand{\E}{\mathbb{E}}
\newcommand{\h}{\mathcal{H}}
\newcommand{\pr}{\mathbb{P}}
\newcommand{\1}{\mathbf{1}}
\newcommand{\da}{\downarrow}
\newcommand{\D}{\mathrm{d}}

\makeatother

\begin{document}
\begin{frontmatter}

\title{On the local approximation of mean densities of random closed sets}
\runtitle{Mean densities of random sets}

\begin{aug}
\author{\fnms{Elena} \snm{Villa}\corref{}\ead[label=e1]{elena.villa@unimi.it}}
\runauthor{E. Villa} 
\address{Department of Mathematics, Universit\`a degli Studi di Milano, via
Saldini 50, 20133 Milano, Italy. \printead{e1}}
\end{aug}

\received{\smonth{2} \syear{2012}}
\revised{\smonth{6} \syear{2012}}

%
\begin{abstract}
Mean density of lower dimensional random closed sets, as well as
the mean boundary density of
full dimensional random sets, and their
estimation are of great interest in many real applications. Only
partial results are available so far in current literature,
under the assumption that the random set is either stationary, or it is
a Boolean model, or it has convex grains. We consider here
non-stationary random closed sets (not necessarily Boolean
models), whose grains have to satisfy some general regularity
conditions, extending previous results. We address the open problem
posed in
(\textit{Bernoulli} \textbf{15} (2009) 1222--1242)
about the approximation of the mean density of lower dimensional
random sets by a pointwise
limit, and to the open problem posed by Matheron in
(\textit{Random Sets and Integral Geometry} (1975) Wiley)
about the existence (and its value) of the
so-called specific area of full dimensional random closed sets.
The relationship with the spherical contact distribution
function, as well as
some examples and applications are also discussed.
\end{abstract}

%
\begin{keyword}
\kwd{mean density}
\kwd{Minkowski content}
\kwd{random measure}
\kwd{specific area}
\kwd{stochastic geometry}
\end{keyword}

\end{frontmatter}

\section{Introduction}\label{secintro} We remind that a \textit{random
closed set} $\Theta$ in
$\mathbb{R}^d$ is a measurable map
\[
\Theta\dvtx (\Omega,\mathfrak{F},\mathbb{P})\longrightarrow(\mathbb {F},
\sigma_\mathbb{F}),
\]
where $\mathbb{F}$ denotes the class of the closed subsets in
$\mathbb{R}^d$, and $\sigma_\mathbb{F}$ is the $\sigma$-algebra
generated by the so called \emph{Fell topology}, or
\emph{hit-or-miss topology}, that is the topology generated by the
set system
\[
\{\mathcal{ F}_G\dvt G\in\mathcal{ G}\} \cup \bigl\{\mathcal{
F}^C\dvt C\in\mathcal{ C}\bigr\},
\]
where $\mathcal{ G}$ and $\mathcal{ C}$ are the system
of the open and compact subsets of $\R^d$, respectively (e.g., see~\cite{matheron}).
We say that
a random closed set
$\Theta\dvtx (\Omega,\mathfrak{F},\pr)\to(\mathbb{F},\sigma_\mathbb{F})$
satisfies a certain property (e.g., $\Theta$ has
Hausdorff dimension $n$) if $\Theta$ satisfies that property
$\pr$-a.s.; throughout the paper we shall deal with countably
$\h^n$-rectifiable random closed sets. For a discussion about
measurability of $ \mathcal{H}^n(\Theta)$, we refer to
\cite{Zahle,MB}.

Let $\Theta_n$ be a set of locally finite $\h^n$-measure; then it
induces a random measure $\mu_{\Theta_n}$ defined by
\[
\mu_{\Theta_n}(A):=\mathcal{H}^n(\Theta_n\cap A), \qquad A
\in\mathcal{B}_{\mathbb{R}^d},
\]
and the corresponding
expected measure
\[
\mathbb{E}[\mu_{\Theta_n}](A):=\mathbb{E}\bigl[\mathcal{H}^n(
\Theta_n\cap A)\bigr],\qquad  A\in\mathcal{B}_{\mathbb{R}^d}.
\]
Whenever
$\mathbb{E}[\mu_{\Theta_n}]$ is
absolutely continuous with respect to $\h^d$, its density (or
Radon--Nikodym derivative)
with respect to $\h^d$ is called \textit{mean
density of} $\Theta_n$, and it is denoted by
$\lambda_{\Theta_n}$.

The problem of the evaluation and the estimation of the mean
density of lower dimensional random closed sets (i.e., with
Hausdorff dimension less than $d$), and in particular of the mean
surface density $\lambda_{\partial\Theta}$ for full dimensional
random sets, is of great interest in several real applications.
We mention, for instance, applications in image
analysis (e.g., \cite{Galerne} and reference therein),
in medicine (e.g., in studying tumor growth \cite{anderson}),
and in material science in phase-transition models (e.g.,
\cite{ECS10}). (See also \cite{AKV,Benes-Rataj,KV-IAS} and
references therein.)

In particular, we recall that in the
well-known seminal book by Matheron on random closed sets
\cite{matheron}, page~50, the so-called \textit{specific area}
$\sigma_\Theta$ is defined by
%
\begin{equation}
\sigma_{\Theta}(x):=\lim_{r\da
0}\frac{\pr(x\in\Theta_{\oplus r}\setminus\Theta)}{r} ,\label{specific}
\end{equation}
where $\Theta_{\oplus r}$ is the parallel set of $\Theta$ at distance
$r>0$, that is, $\Theta_{\oplus r}:=\{s\in\R^d \dvt \operatorname{ dist}(x,\Theta)\leq r\}$;
it
is introduced as
a probabilistic version of the derivative at
0 of the volume function $V(r):=\h^d(\Theta_{\oplus r})$, and so,
whenever the limit exists, as a possible approximation of what we
denote by $\lambda_{\partial\Theta}$, the mean boundary density
of $\Theta$. The problem of the existence of $\sigma_\Theta$ is
left as an open problem in~\cite{matheron} (apart from particular
cases as stationary random
closed sets).

More recently, in \cite{AKV} the problem of the approximation of
the mean density $\lambda_{\Theta_n}$ of lower dimensional
non-stationary random closed sets is faced under quite general
regularity assumptions on the rectifiability of $\Theta_n$. More
precisely, an approximation of $\lambda_{\Theta_n}$ in weak form
is proved in~\cite{AKV}, Theorem~4; namely
%
\begin{equation}
\lim_{r\downarrow
0}\int_A \frac{\mathbb{P}(x\in\Theta_{n_{\oplus r}})}{b_{d-n}r^{d-n}} \,\mathrm{ d} x=
\int_A \lambda_{\Theta_n}(x) \,\mathrm{ d} x.
\label{scambiostar}
\end{equation}
The possibility of exchanging limit
and integral in the above expression when $\Theta_n$ is not
stationary with $n>0$, was left as open problem in
\cite{AKV}, Remark 8. (The stationary and the 0-dimensional cases
are trivial.)

A first attempt to
solve the above mentioned open problems (the one for
$\sigma_\Theta$ posed by Matheron, and the one for
$\lambda_{\Theta_n}$ with $n<d$ posed in \cite{AKV}), is given in
\cite{MeandensityBool}, where explicit
results are proven for inhomogeneous Boolean models.

The aim of the present paper is to address such open problems for
more general random closed sets. Indeed, even if Boolean models
are widely studied in stochastic geometry (e.g., see \cite{CIME}),
it is clear that they cannot be taken as model for many real
situations in applications. Thus, we revisit here some results in
\cite{MeandensityBool}, addressing the two mentioned open
problems; we provide sufficient conditions on lower dimensional
random sets $\Theta_n$ so that
%
\begin{equation}
\label{aim} \lambda_{\Theta_n}(x)=\lim_{r\downarrow0}\frac{\mathbb{P}(x\in
{\Theta_{n_{\oplus r}}})}{b_{d-n}r^{d-n}},
\qquad \h^d\mbox{-a.e. }x\in\R^d,
\end{equation}
and so that the specific area $\sigma_\Theta$ defined as limit in
\eqref{specific} exists, in the case of random sets $\Theta$ with
non-negligible $\h^d$-measure. Such results might allow to face a
wider class of possible applications; indeed, for instance, the
statistical estimator $\widehat{\lambda}_{\Theta_n}^N(x)$ of the
mean density $\lambda_{\Theta_n}(x)$, introduced in \cite{MeandensityBool} and which we recall here in Corollary \ref
{corstimatore}, can now be applied to very general lower dimensional
random sets $\Theta_n$, not only in stationary settings or to
Boolean models, and so also to non-stationary germ-grains model
whose grains are not assumed to be independent.
We also mention here that the estimation of $\lambda_{\Theta_n}$
and $\sigma_{\Theta}$ might be considered as the stochastic
analogous to the estimation of a non-random unknown support, and the
stochastic counterpart of boundary estimation for a given
support, respectively (see, e.g., \cite{cuevas1,cuevas2}); this
might lead to possible further research on this topics.

The plan of the paper is the following: preliminary notions and
known results on the so-called Minkowski content of sets and on
point processes and germ-grain models are briefly recalled in
Section \ref{sec2}. In Section \ref{sec3}, we answer to the open problem posed in
\cite{AKV} mentioned above, that is we prove equation \eqref{aim}; we
also provide an explicit expression for $\lambda_{\Theta_n}(x)$. A
natural estimator follows as a corollary. Further results and
remarks are discussed in the final part of the section; known
results on the special case of Boolean models follow here as
particular case. In Section~\ref{sec4}, random sets with non-negligible
$\h^d$-measure are considered; by recalling recent results on the
outer Minkowski content notion we answer to the open problem posed
by Matheron in \cite{matheron} about the existence of the specific
area $\sigma_\Theta$ of random sets $\Theta$ which can be
represented as \emph{one-grain} random sets. The relationship
between $\sigma_\Theta$, the mean boundary density
$\lambda_{\partial\Theta}$ of $\Theta$, and its spherical contact
distribution function is studied. Some explicit formulas for the
derivative of the contact distribution are also proved.

\section{Preliminaries and notation}\label{sec2}
In this section, we recall basic definitions, notation and results
on point processes and geometric measure theory which we shall use
in the following.

\subsection{The Minkowski content notion and related results}
Throughout the paper, $\h^n$ is the $n$-dimensional Hausdorff
measure, $\D x$ stands for $\h^d(\D x)$, and~$\mathcal{B}_\mathcal{X}$ is the Borel $\sigma$-algebra of any
space $\mathcal{X}$.
$B_r(x)$, $b_n$ and $\mathbf{S}^{d-1}$
will denote the closed ball with\vadjust{\goodbreak} centre $x$ and radius $r\geq0$,
the volume of the unit ball in $\mathbb{R}^n$ and the unit sphere
in $\R^d$, respectively. We remind that a compact set
$A\subset\mathbb{R}^d$ is called \emph{$n$-rectifiable} ($0\leq
n\leq d-1$ integer) if it can be written as the image of a
compact subset of $\mathbb{R}^n$ by a Lipschitz map from $\R^n$ to
$\R^d$; more in general, a closed subset $A$ of $\mathbb{R}^d$ is
said to be \emph{countably $\mathcal{H}^{n}$-rectifiable} if there
exist countably many $n$-dimensional Lipschitz graphs
$\Gamma_i\subset\mathbb{R}^d$ such that $A\setminus\bigcup_i\Gamma_i$
is $\mathcal{H}^n$-negligible. (For definitions and basic
properties of Hausdorff measure and rectifiable sets see, e.g.,
\cite{AFP,Falconer,Fed-libro}.)

The notion of \textit{$n$-dimensional Minkowski content} will play a
fundamental role
throughout the paper. We recall that, given a subset $A$ of $\R^d$
and an integer $n$ with $0\leq n\leq d$, the \emph{$n$-dimensional
Minkowski content of $A$} is defined as
%
\begin{equation}
\mathcal{M}^n(A):=\lim_{r\da0}\frac{\h^d(A_{\oplus
r})}{b_{d-n}r^{d-n}},
\label{minkcont}
\end{equation}
whenever the limit
exists finite. Well known general results about the existence of
the Minkowski content of closed sets in $\mathbb{R}^d$ are
related to rectifiability properties of the involved
sets. In particular, the following theorem is proved in
\cite{AFP}, page~110. (We call Radon measure in $\mathbb{R}^d$
any non-negative and $\sigma$-additive set function
defined on
$\mathcal{B}_{\mathbb{R}^d}$ which is finite on bounded
sets.)
%
\begin{theorem}\label{teominkcont} Let
$A\subset\mathbb{R}^d$ be a countably
$\mathcal{H}^{n}$-rectifiable compact set, and assume that
%
\begin{equation}
\eta\bigl(B_r(x)\bigr)\geq\gamma r^n \qquad\forall x\in A,\
\forall r\in(0,1) \label{condizioneeta}
\end{equation}
holds for some $\gamma>0$ and
some Radon measure $\eta\ll\h^n$ in $\mathbb{R}^d$. Then
$\mathcal{M}^n(A)=\mathcal{H}^n(A)$.
\end{theorem}
Condition
\eqref{condizioneeta} is a kind of quantitative non-degeneracy
condition which prevents $A$ from being too sparse; simple
examples show that $\mathcal{M}^n(A)$ can be infinite, and
$\mathcal{H}^{n}(A)$ arbitrarily small, when this condition fails
\cite{AFP,AColeV}. The above theorem extends (see \cite{AFP}, Theorem 2.106) the
well-known Federer's result \cite{Fed-libro}, page 275, to
countably $\mathcal{H}^{n}$-rectifiable compact sets; in
particular for any $n$-rectifiable compact set $A\subset\R^d$
there exists a suitable measure $\eta$ satisfying
\eqref{condizioneeta} (see \cite{AColeV}, Remark~1). As a
consequence, for instance in the case $n=d-1$, the boundary of
any convex body or, more in general, of a set with positive
reach, and the boundary of a set with Lipschitz boundary satisfy
condition \eqref{condizioneeta}. Note also that if a Radon
measure $\eta$ as in Theorem \ref{teominkcont} exists, then
it can be assumed to be a probability
measure, without loss of generality (e.g., see \cite{MeandensityBool}); the next theorem is proved in \cite{MeandensityBool},
and provides a result on the existence of the limit in \eqref
{minkcont} when the measure $\h^d$ is replaced by a measure having
density $f$ with respect to $\h^d$, and so it may be seen as a
generalization of the theorem above. $\operatorname{ disc}f$ denotes the set
of all the points of discontinuity of $f$.
%
\begin{theorem}\label{Minkgen} Let $\mu\ll\h^d$ be a positive
measure in $\R^d$, admitting a
locally bounded density $f$, and $A\subset\mathbb{R}^d$ be a
countably $\mathcal{H}^{n}$-rectifiable compact set such that
condition \eqref{condizioneeta} holds for some $\gamma>0$ and
some probability measure $\eta\ll\h^n$ in $\R^d$. If $\h^n(\operatorname{
disc}f)=0$, then
\[
\lim_{r\da0}\frac{\mu(A_{\oplus r})}{b_{d-n}r^{d-n}}=\int_A f(x)
\h^n(\D x).
\]
\end{theorem}

\subsection{Point processes}
Here we report some known facts from the theory of point processes
just for establishing notation which will be used later. For a
more complete exposition of the theory of point processes, see,
for example, \cite{daley-vj}. Roughly speaking a point process
$\widetilde{\Phi}$ in $\R^d$ is a locally finite collection
$\{\xi_i\}_{i\in\N}$ of random points in $\R^d$. Formally,
$\widetilde{\Phi}$ can be seen as a random counting measure, that
is a measurable map from a probability space
$(\Omega,\mathcal{F},\pr)$ into the space of locally finite
counting measures on $\R^d$. $\widetilde{\Phi}$ is called simple
if $\widetilde{\Phi}(\{x\})\leq1$ for all $x\in\R^d$; we shall
always consider simple point processes.

The measure $\widetilde{\Lambda}(A):=\E[\widetilde{\Phi}(A)]$ on
$\mathcal{B}_{\R^d}$ is called \textit{intensity measure} of $\widetilde{\Phi}$; whenever it is
absolutely continuous with respect to $\h^d$, its density is
called \textit{intensity} of $\widetilde{\Phi}$. It is well known the
so-called Campbell's formula (e.g., see \cite{CIME}, page~28), which
states that for any measurable function $f\dvtx \R^d\to\R$ the
following holds
\[
\E \biggl[\sum_{x\in\widetilde{\Phi}}f(x) \biggr]=\int
_{\R
^d}f(x)\widetilde {\Lambda}(\D x).
\]
Another important measure associated to a point process
$\widetilde{\Phi}$ is the so-called \textit{second factorial moment
measure} $\widetilde{\nu}_{[2]}$ of $\widetilde{\Phi}$; it is the
measure on $\mathcal{B}_{\R^{2d}}$ defined by (e.g., see
\cite{CIME,SKM})
\[
\int f(x,y)\widetilde{\nu}_{[2]}\bigl(\D(x,y)\bigr)=\E \biggl[\sum
_{x,y\in
\widetilde
{\Phi},   x\neq y}f(x,y) \biggr]
\]
for any non-negative
measurable function $f$ on $\R^{2d}$. Moreover, $\widetilde{\Phi}$
is said to have \textit{second moment density} $\widetilde g$ if
$\widetilde{\nu}_{[2]}=\widetilde{g}\nu^{2d}$, that is
\[
\widetilde{\nu}_{[2]}(C)=\int_C
\widetilde{g}(x,y)\,\D x\,\D y
\]
for any compact $C\subset\R^{2d}.$ Informally,
$\widetilde{g}(x,y)$ represents the joint probability that there
are points at two specific locations $x$ and $y$:
\[
\widetilde{g}(x,y)\,\D x\,\D y\sim\pr\bigl(\widetilde{\Phi}(\D x)>0, \widetilde{
\Phi}(\D y)>0\bigr).
\]

A generalization of the above notion is the so-called \textit{marked
point process}. We recall that a marked point process
$\Phi=\{\xi_i,K_i\}_{i\in\N}$ on $\R^d$ with marks in a complete
separable metric space (c.s.m.s.)
$\mathbf{K}$ is a point
process on $\R^d\times\mathbf{K}$ with the property that the unmarked
process $\{\widetilde{\Phi}(B)\dvt B\in\mathcal{B}_{\R^d} \}:=\{\Phi
(B\times\mathbf{K})\dvt B\in\mathcal{B}_{\R^d}\}$ is a point
process in $\R^d$. $\mathbf{K}$ is called \emph{mark space}, while
the random element $K_i$ of $\mathbf{K}$ is the \emph{mark
associated to the point} $\xi_i$.
$\Phi$ is said to be \textit{stationary} if the distribution of $\{\xi_i+x,K_i\}_i$ is
independent of $x\in\R^d$.

If the marks are independent and identically distributed, and
independent of the unmarked point process $\widetilde{\Phi}$, then
$\Phi$ is said to be an \textit{independent marking} of
$\widetilde{\Phi}$.

The intensity measure of $\Phi$, say $\Lambda$, is a $\sigma$-finite
measure on $\mathcal{B}_{
\R^d\times\mathbf{K}}$ defined as $\Lambda(B\times
L):=\mathbb{E}[\Phi(B\times L)]$, the mean number of points of
$\Phi$ in $B$
with marks in $L$.
We recall that a Campbell's formula for marked point processes
holds as well \cite{CIME}:
%
\begin{equation}
\E \biggl[\sum_{(x,K)\in\Phi}f(x,K) \biggr]=\int
_{\R^d\times\mathbf
{K}}f(x,K)\Lambda\bigl(\D(x,K)\bigr).\label{Campformula}
\end{equation}
Since $\mathbf
{K}$ is a c.s.m.s. and
$\widetilde{\Lambda}$ is a $\sigma$-finite measure, it is possible
to factorize $\Lambda$ in the following way \cite{Last-Brandt}:
\[
\Lambda\bigl(\D(x,K)\bigr)=\kappa(x,\D K)\widetilde{\Lambda}(\D x),
\]
where $\widetilde{\Lambda}$
is the intensity measure of
the unmarked process $\widetilde{\Phi}$, and $\kappa(x, \cdot)$
is a probability measure on $\mathbf{K}$ for all $x\in\R^d$,
called the mark distribution at point $x$. A common assumption
(e.g., see \cite{Hug-Last}) is that there exist a measurable
function $\lambda\dvtx \R^d\times\mathbf{K}\to\R_+$ and a probability
measure $Q$ on $\mathbf{K}$ such that
%
\begin{equation}
\Lambda\bigl(\D(x,K)\bigr)=\lambda(x,K)\,\D x Q(\D K), \label{star}
\end{equation}
this
happens if and only if $\kappa(x,\cdot)$ is absolutely
continuous with respect to $Q$ for $\h^d$-a.e. $x\in\R^d$.

If $\Phi$ is stationary, then its intensity measure is of the type
$\Lambda=\lambda\nu^d\otimes Q$ for some $\lambda>0$ and $Q$
probability measure on $\mathbf{K}$. If $\Phi$ is an independent
marking of $\widetilde{\Phi}$, then
$\Lambda(\D(x,K))=\widetilde{\Lambda}(\D x)Q(\D K) $, where $Q$ is
a probability measure on $\mathbf{K}$, called
distribution of the marks.

Let $(\R^{d}\times\mathbf{K})^2:=\R^{d}\times\mathbf{K}\times
\R^{d}\times\mathbf{K}$; the \textit{second factorial moment measure}
$\nu_{[2]}$ of $\Phi$ is the measure on
$\mathcal{B}_{(\R^{d}\times\mathbf{K})^2}$ so defined
\cite{SKM}
%
\begin{equation}
\int f(x_1,K_1,x_2,K_2)
\nu_{[2]}\bigl(\D(x_1,K_1,x_2K_2)
\bigr)= \E \biggl[\mathop{\sum_{(x_i,K_i),(x_j,K_j)\in\Phi,}}_{
x_i\neq x_j}f(x_i,K_i,x_j,K_j)
\biggr]\label{nu2def}
\end{equation}
for any non-negative
measurable function $f$ on
$(\R^{d}\times\mathbf{K})^2$. By denoting $\widetilde{\nu}_{[2]}$
the second factorial moment measure of the unmarked process
$\widetilde{\Phi}$, for any $B_1,B_2\in\mathbf{K}$ the measure
$\nu_{[2]}(\cdot\times B_1\times\cdot\times B_2)$ is absolutely
continuous with respect to $\widetilde{\nu}_{[2]}$; moreover, if
$\widetilde{\nu}_{[2]}$ is $\sigma$-finite then
%
\begin{equation}
{\nu}_{[2]}\bigl(\D(x_1,K_1,x_2,K_2)
\bigr)=M_{x_1,x_2}\bigl(\D(K_1,K_2)\bigr)\widetilde {
\nu }_{[2]}\bigl(\D(x_1,x_2)\bigr),
\label{nu2gen}
\end{equation}
where $M_{x_1,x_2}$ is a measure on
$\mathbf{K}^2$ for any fixed $x_1$ and $x_2$, called \textit{two-point mark distribution}. Informally,
${\nu}_{[2]}(\D(x_1,K_1,x_2,K_2))$ represents the joint
probability that there are points at two specific locations
$x_1$ and $x_2$ with marks $K_1$ and $K_2$, respectively.

Similarly to $\Lambda$, we shall assume that there exist a
measurable function $g\dvtx (\R^d\times\mathbf{K})^2\to\R_+$ and a
probability measure $Q_{[2]}$ on $\mathbf{K}^2$ such that
%
\begin{equation}
\nu_{[2]}\bigl(\D(x_1,K_1,x_2,K_2)
\bigr)=g(x_1,K_1,x_2,K_2)\,\D
x_1\,\D x_2 Q_{[2]}\bigl(\D(K_1,K_2)
\bigr).\label{nu2}
\end{equation}
We remind that if $\Phi$ is a
marked Poisson point process with intensity measure
$\Lambda(\D(x,K))=\kappa(x,\D K)\widetilde{\Lambda}(\D x)$, then $
\widetilde{\nu}_{[2]}
=\widetilde{\Lambda}\otimes\widetilde{\Lambda}$ and
$\nu_{[2]}=\Lambda\otimes\Lambda$, and so
\[
M_{x,y}\bigl(\D(s,t)\bigr)=\kappa(x,\D s)\kappa(y,\D t) ;
\]
in particular, by the assumptions \eqref{star} and
\eqref{nu2} it follows
%
\begin{eqnarray}\label{poisson}
g(x_1,K_1,x_2,K_2)&=&
\lambda(x_1,K_1)\lambda(x_2,K_2),
\nonumber
\\[-8pt]
\\[-8pt]
\nonumber
Q_{[2]}\bigl(\D(K_1,K_1)\bigr)&=&Q(\D
K_1)Q(\D K_2).
\end{eqnarray}

We also recall that point processes can be considered on quite
general metric spaces. In particular, a point process in
$\mathcal{C}^d$, the class of compact subsets of $\R^d$, is called
\textit{particle process} (e.g., see \cite{CIME} and references
therein). It is well known that, by a \textit{center map}, a particle
process can be transformed into a marked point process $\Phi$ on
$\R^d$ with marks in $\mathcal{C}^d$, by representing any compact
set $C$ as a pair $(x,Z)$, where $x$ may be interpreted as the
``location'' of $C$ and $Z:=C-x$ the ``shape'' (or ``form'') of
$C$ (e.g., see \cite{CIME}, page~192 and \cite{HLW-survey}). In this
case the marked point process $\Phi=\{(X_i,Z_i)\}$ is also called
\textit{germ-grain model}. In case of independent marking, the grains
$Z_i$'s are i.i.d. as a typical grain $Z_0$ with mark
distribution $Q$, which is also called, in this case, \textit{grain
distribution} or \textit{distribution of the typical grain}.

Every random closed set in $\R^d$ can be represented as a
germ-grain model, and so by a suitable marked point process
$\Phi=\{X_i,Z_i\}$. In many examples and applications the random
sets $Z_i$ are uniquely determined by suitable random parameters
$S\in\mathbf{K}$. For instance, in the very simple case of random
balls, $\mathbf{K}=\R_+$ and $S$ is the radius of a ball centred
in the origin;
in applications to birth-and-growth processes, in some
models $\mathbf{K}=\R^d$ and $S$ is the spatial location of the
nucleus (e.g., \cite{AKV}, Example 2); in segment processes in
$\R^2$, $\mathbf{K}=\R_+\times[0,2\pi]$ and $S=(L,\alpha)$ where~$L$ and $\alpha$ are the random length and orientation of the
segment through the origin, respectively (e.g.,~\cite{MeandensityBool}, Example 2); etc. So, in order to use
similar notation to previous works (e.g.,~\cite{MeandensityBool,ECS10}), we shall consider random sets $\Theta$
described by
marked point processes $\Phi=\{(X_i,S_i)\}$ in $\R^d$ with marks
in a suitable mark space $\mathbf{K}$ so that $Z_i=Z(S_i)$ is a
random set containing the origin:
%
\begin{equation}
\Theta(\omega)=\bigcup_{(x_i,s_i)\in\Phi(\omega
)}x_i+Z(s_i),\qquad
\omega\in\Omega .\label{theta}
\end{equation}
We also recall that whenever $\Phi$ is a marked
Poisson point process, $\Theta$ is said to be a \textit{Boolean
model}.

The intensity measure $\Lambda$ of $\Phi$ is
commonly assumed to be such that the mean number of grains hitting
any compact subset of $\R^d$ is finite, which is equivalent to say
that the mean number of grains hitting the ball $B_R(0)$ is finite
for any $R>0$:
%
\begin{equation}
\E\biggl[\sum_{(x_i,s_i)\in\Phi}\1_{(-Z(s_i))_{\oplus
R}}(x_i)
\biggr]\stackrel{\eqref{Campformula}} {=}\int_{\R^d\times\mathbf{K}}
\1_{(-Z(s))_{\oplus R}}(x)\Lambda\bigl(\D(x,s)\bigr) <\infty \qquad \forall R>0.\label{cond}
\end{equation}
%

\section{Mean densities of lower dimensional random closed sets}\label{sec3}

\subsection{Assumptions}\label{secass}

Let $\Theta_n$ be a random closed set in $\R^d$ with integer
Hausdorff dimension $0<n<d$ as in \eqref{theta}, where
$\Phi$ has intensity measure $\Lambda(\D(x,s))=\lambda(x,s)\,\D
xQ(\D s)$ and second factorial moment measure
$\nu_{[2]}(\D(x,s,y,t))=g(x,s,y,t)\,\D x \,\D y Q_{[2]}(\D(s,t))$ such
that the following assumptions are fulfilled:
\begin{longlist}[(A2)]
\item[(A1)] for any $(y,s)\in\R^d\times\mathbf{K}$, $y+Z(s)$ is a
countably $\h^n$-rectifiable and compact subset of $\R^d$, such
that there exists a closed set $\Xi(s)\supseteq Z(s)$ such that
$\int_\mathbf{K}\h^n(\Xi(s))Q(\D s)<\infty$ and
%
\begin{equation}
\h^n\bigl(\Xi(s)\cap B_r(x)\bigr)\geq\gamma
r^n\qquad  \forall x\in Z(s),\forall r\in(0,1) \label{gamma}
\end{equation}
for some $\gamma>0$ independent of $y$ and $s$;

\item[(A2)] for any
$s\in\mathbf{K} $, $\h^n(\operatorname{ disc}(\lambda(\cdot,s)))=0$ and
$\lambda
(\cdot,s)$ is locally bounded
such that
for any compact $K\subset\R^d$
%
\begin{equation}
\sup_{x\in K_{\oplus\operatorname{
diam}(Z(s))}}\lambda(x,s)\leq\widetilde{\xi}_K(s)\label{locbound1}
\end{equation}
for some $\widetilde{\xi}_K(s)$ with
$\int_{\mathbf{K}}\h^n(\Xi(s))\widetilde{\xi}_K(s)Q(\D s)<\infty$;

\item[(A3)] for any
$(s,y,t)\in\mathbf{K}\times\R^d\times\mathbf{K}$, $\h^n(\operatorname{
disc}(g(\cdot,s,y,t)))=0$ and $g(\cdot,s,y,t)$ is locally bounded
such that
for any compact $K\subset\R^d$ and $a\in\R^d$,
%
\begin{equation}
\1_{(a-Z(t))_{\oplus1}}(y)\sup_{x\in{K}_{\oplus\operatorname{
diam}(Z(s))}}g(x,s,y,t)\leq\xi_{a,K}(s,y,t)
\label{locbound}
\end{equation}
for some $\xi_{a,K}(s,y,t)$ with
$\int_{\R^d\times\mathbf{K}^2}\h^n(\Xi(s))\xi_{a,K}(s,y,t)\,\D
yQ_{[2]}(\D s,\D t)<\infty$.
\end{longlist}

Before stating our main results, we briefly
discuss the above assumptions. As mentioned in the \hyperref[secintro]{Introduction},
we want to find sufficient conditions such that
equation \eqref{aim} holds for a general class of random closed sets
$\Theta_n$, so answering to the open problem stated in
\cite{AKV}, Remark 8. We point out that such a result has been
proved recently in \cite{MeandensityBool} for Boolean models
with position-independent grains, and so only in the case in which
$\Phi$ is a Poisson point process with intensity measure $\Lambda$
of the type $\Lambda(\D((x,s)))=\lambda(x)\,\D xQ(\D s)$. In that
work, the assumption that $\Phi$ was a marked Poisson point
process allowed to apply the explicit expression of the capacity
functional of $\Theta_n$, both in proving the exchange between limit
and integral in \eqref{scambiostar}, and in providing an explicit
formula for the mean density $\lambda_{\Theta_n}$ of $\Theta_n$ in
terms of its intensity measure $\Lambda$. Actually, in order to
prove equation \eqref{aim}, the knowledge of the capacity functional of
$\Theta_n$ is not necessary, by making use of Campbell's formula.
Nevertheless, for a general random set $\Theta_n$ as in the above
assumptions, and so without the further assumption that $\Phi$ is
a marked Poisson process, we need to introduce also the second
factorial moment measure of $\Phi$, and the related assumption
(A3). Of course, considering here a generic random set $\Theta_n$
(point process $\Phi$), it obvious that the above assumptions\vadjust{\goodbreak}
are similar to (actually, they generalize) those which appear in
\cite{MeandensityBool}; as a matter of fact (A1$'$) and (A2$'$) in
\cite{MeandensityBool} coincide with (A1) and (A2) above in the
case of independent marking. We also point out that in the
particular case of Boolean models, the second factorial moment
measure $\nu_{[2]}$ is given in terms of the intensity measure
$\Lambda$, and so the function $g$ in terms of $\lambda$ by
\eqref{poisson}; this is the reason why here assumption (A3)
appears, whereas it is already contained in (A1$'$) and (A2$'$) in
\cite{MeandensityBool}, Theorem 3.13 (see also Corollary \ref{corBoolean}
below).

We mention also that taking $ \nu_{[2]}$ of
the type $\nu_{[2]}(\D(x,s,y,t))\,{=}\,g(x,s,y,t)\,\D x \,\D y
Q_{[2]}(\D(s,t))$ is in accordance to the assumption in
\cite{Hug-Last}, Proposition~4.9, where contact distributions of general
germ-grain models with compact convex grains are considered; in
that paper $\nu_{[2]}$ is assumed to be absolutely continuous with
respect to the product measure $\h^d\otimes\mu$, where $\mu$ is
$\sigma$-finite measure on $\mathbf{K}\times\R^d\times\mathbf{K}$.

Moreover,
note that the measure $\h^n(\Xi(s)\cap \cdot)$ in (A1) plays the
same role as the measure~$\eta$ of Theorem \ref{teominkcont};
indeed (A1) might be seen as the stochastic version of
\eqref{condizioneeta}. (See also \cite{MeandensityBool}, Remark~3.6,
and the examples discussed in \cite{AKV}.) Roughly
speaking, such an assumption tells us that each possible grain
associated to any point $x$ of the underling point process
$\widetilde{\Phi}$ is sufficiently regular, so that it admits
$n$-dimensional Minkowski content; this explains also why
requiring the existence of a constant $\gamma$ as in (A1)
independent on $y$ and $s$ is not too restrictive (see also the
example below about this). Note that the condition
$\int_\mathbf{K}\h^n(\Xi(s))Q(\D s)<\infty$ means that the
$\h^n$-measure of the grains is finite in mean. In order to
clarify better the meaning of assumption (A1), let us consider
the following simple example.
%
\begin{example}\label{exsegment}
Let $\Theta_1$ be a germ grain model with segments as grains,
with random length. (As it will be clear, the orientation of
the segments does not take part to the validity of (A1).)
Let us only assume that the mean length of the grain is
finite.
We may notice that the introduction of the suitable random set
$\Xi$ is needed only if the length of the segments could be
indefinitely close to~0. Indeed, let us first consider the case
in which the length is bounded from below by a positive constant, for instance
$\h^1(Z(s))\geq l>0$ for any $s\in\mathbf{K}$; then
\[
\h^1\bigl(Z(s)\cap B_r(x)\bigr)\geq\min\{l, 1\}r\qquad
\forall x\in Z(s),\  \forall r\in(0,1),
\]
and so there exists $\gamma:=\min\{l, 1\}>0$, clearly independent
of the position and of the length of the particular grain
considered.

Now let us consider the case in which the length is not bounded
from below by a positive constant (e.g., the length is uniformly
distributed in $[0,L]$). In this case, $l=0$ and so we have to
introduce a suitable random set $\Xi$ satisfying \eqref{gamma}; a
possible solution is to extend all the segments having length
less than 2 (the extension can be done homothetically from the
center of the segment, so that measurability of the process is
preserved). In particular, for any $s\in\mathbf{K}$, let
\[
\Xi(s)=\cases{
Z(s), & \quad $\mbox{if }
\mathcal{H}^1\bigl(Z(s)\bigr)\geq2,$
\vspace*{2pt}\cr
Z(s) \mbox{ extended to length }2,&\quad  $\mbox{if }\mathcal{H}^1
\bigl(Z(s)\bigr)<2;$ }
\]
it follows that \eqref{gamma} holds now with $\gamma=1$.
Since we have assumed that the mean length of the segments is finite,
it follows that $\int_\mathbf{K}\h^n(\Xi(s))Q(\D s)<\infty$,
and so (A1) is fulfilled.

Note that we have chosen segments as grains in order to make
the example simpler, but it is now clear that the same argument
may applied to fibre processes (in order to provide another example of
a random closed set of dimension 1), or even more complicated
random sets in $\R^d$ with any integer dimension $n$.
\end{example}
The role of assumption (A2) and (A3) is more technical, and it
will be clearer later in the proofs of the next statements.
Finally, it is clear that if $\lambda$ and $g$ are bounded, the
above assumptions (A2) and (A3) simplify (see also Remark \ref{indmark}).

\subsection{Main theorem and related results}
In this section, we state and prove our main theorem
(Theorem \ref{thlambda-lim}), which provides a pointwise limit
representation of the mean density $\lambda_{\Theta_n}$ of
$\Theta_n$.
To this aim we need to prove some other related results,
before. We start with the following lemma, which tells us that the
grains of the random set $\Theta_n$ overlap only on a set having
negligible $\h^n$-measure in mean.

\begin{lemma}\label{lemmaattesasomma} Let $\Theta_n$ be a random
closed set in $\R^d$ with
integer
Hausdorff dimension $0<n<d$ as in \eqref{theta}, where
$\Phi$ has intensity measure $\Lambda(\D(x,s))=\lambda(x,s)\,\D
xQ(\D s)$ and second factorial moment measure
$\nu_{[2]}(\D(x,s,y,t))=g(x,s,y,t)\,\D x \,\D y Q_{[2]}(\D(s,t))$.
Then
\[
\E \biggl[\mathop{\sum_{(y_i,s_i),(y_j,s_j)\in\Phi,}}_{
y_i\neq y_j} \h^n\bigl(
\bigl(y_i+Z(s_i)\bigr)\cap\bigl(y_j+Z(s_j)
\bigr)\bigr) \biggr]=0.
\]
\end{lemma}
\begin{pf}
The following chain of equalities hold:
\begin{eqnarray*}
&&\E \biggl[\mathop{\sum_{(y_i,s_i),(y_j,s_j)\in\Phi,}}_{
y_i\neq y_j
} \h^n\bigl(
\bigl(y_i+Z(s_i)\bigr)\cap\bigl(y_j+Z(s_j)
\bigr)\bigr) \biggr]
\\
&&\quad\stackrel{ \eqref {nu2def}} {=}\int_{(\R^d\times\mathbf{K})^2}\h^n
\bigl(\bigl(x+Z(s)\bigr)\cap \bigl(y+Z(t)\bigr)\bigr)\nu_{[2]}\bigl(\D
(x, s, y, t)\bigr)
\\
&&\quad=\int_{(\R^d\times\mathbf{K})^2} \biggl(\int_{\R^d}
\1_{x+Z(s)}(u)\1_{y+Z(t)}(u)\h^n(\D u) \biggr)
\nu_{[2]}\bigl(\D(x, s, y, t)\bigr)
\\
&&\quad= \int_{(\R^d\times\mathbf{K})^2} \biggl(\int_{\R^d}
\1_{u-Z(s)}(x)\1_{u-Z(t)}(y)\h^n(\D u) \biggr)g(x,s,y,t)\,\D
x\,\D y Q_{[2]}(\D s,\D t)
\\
&&\quad=\int_{\R^d} \biggl(\int_\mathbf{K}\int
_{\R^d}\int_{\mathbf{K}}\1_{u-Z(t)}(y) \int
_{u-Z(s)}g(x,s,y,t)\,\D x \,\D yQ_{[2]}(\D s,\D t) \biggr)
\h^n(\D u),
\end{eqnarray*}
where the last equality is implied by Fubini's theorem. The
assertion follows by observing that $\int_{u-Z(s)}g(x,s,y,t)\,\D
x=0$, because $\h^d(Z(s))=0$, being lower dimensional.
\end{pf}

In order to prove our next results, we recall that in \cite{AKV}
it is proved that if $S\subset\R^d$ is a countably
$\h^n$-rectifiable compact set such that
\[
\eta\bigl(B_r(x)\bigr)\geq\gamma r^n \qquad \forall x\in S,\
\forall r\in(0,1)
\]
holds for some $\gamma>0$ and some finite measure $\eta\ll\h^n$
in $\R^d$, then
%
\begin{equation}
\frac{\h^d(S_{\oplus r})}{b_{d-n}r^{d-n}}\leq \frac{\eta(\R^d)}{\gamma}2^n4^d
\frac{b_d}{b_{d-n}}\qquad \forall r<2.\label{condeta1}
\end{equation}

\begin{remark}\label{remradon} By \eqref{condeta1}, and the
proof of
Lemma 3.14 in \cite{MeandensityBool}, we know that
\[
\h^{d}\bigl(Z(s)_{\oplus R}\bigr)\leq\cases{ %
\h^{n}\bigl(\Xi(s)\bigr)\gamma^{-1}2^n4^db_dR^{d-n},
& \quad $\mbox{if }R<2,$
\vspace*{2pt}\cr
\h^{n}\bigl(\Xi(s)\bigr)\gamma^{-1}2^n4^db_dR^{n},
& \quad $\mbox{if }R\geq2,$}
\]
and so condition \eqref{cond}, which guarantees that the mean
number of grains intersecting any compact subset of $\R^d$ is
finite, is fulfilled:
\begin{eqnarray*}
&&\int_{\R^d\times
\mathbf{K}}\1_{(-Z(s))_{\oplus R}}(x)\Lambda\bigl(\D(x,s)\bigr)
\\
&&\quad\leq 2^n4^db_d\max\bigl\{R^{d-n};
R^d\bigr\}\int_{\mathbf{K}}\widetilde{\xi}_{B_R(0)}
\h^n\bigl(\Xi(s)\bigr)Q(\D s)\stackrel{\mathrm{(A2)}} {<}\infty \qquad\forall R>0.
\end{eqnarray*}
As a consequence, together with assumption (A1) which tells us that
each grain has finite
$\h^n$-measure in mean, it is easy to see that
$\E[\mu_{\Theta_n}]$ is locally bounded. Moreover, by proceeding
along the same lines of the proof of Proposition 3.8 in \cite
{MeandensityBool}, we get that $\E[\h^n(\Theta_n\cap A)]=0$ for any
$A\subset\R^d$ with $\h^d(A)=0$, that is $\E[\mu_{\Theta_n}]$ is
absolutely continuous with respect to $\h^d$.
\end{remark}
By following the
hint given in \cite{MeandensityBool}, page~494 (there given for
Boolean models, but here applied to more general $\Theta_n$), the
following proposition, which provides an explicit formula of the
mean density $\lambda_{\Theta_n}$ of $\Theta_n$ in terms of its
intensity measure, is easily proved by means of the above lemma
and Campbell's formula. (See also \cite{lars} for a similar
application.)

\begin{proposition}\label{proplambda} Under the hypotheses of
Lemma \ref{lemmaattesasomma},
%
\begin{equation}
\lambda_{\Theta_n}(y)=\int_{
\mathbf{K}}\int_{y-Z(s)}
\lambda(x,s)\h^n(\D x)Q(\D s)\qquad \mbox{for }\h^d\mbox{-a.e.
}y\in\R^d. \label{formulalambda}
\end{equation}
\end{proposition}

\begin{pf}
By Lemma \ref{lemmaattesasomma}, we know that the event that
different grains of $\Theta_n$ overlap in a subset of $\R^d$ of
positive $\h^n$-measure has null probability; then the following
chain of equalities holds for any $A\in\mathcal{B}_{\R^d}$:
\begin{eqnarray*}
\E\bigl[\h^n(\Theta_n\cap A)\bigr]&=&\E\biggl[\sum
_{(y_i,s_i)\in
\Phi}\h^n\bigl(\bigl(y_i+Z(s_i)
\bigr)\cap A\bigr)\biggr]
\\
&\stackrel{\eqref{Campformula}} {=}&\int_{\R^d\times{
\mathbf{K}}}
\h^n\bigl(y+Z(s)\cap A\bigr)\Lambda\bigl(\D(y,s)\bigr)
\\
&=&\int_{\R^d\times\mathbf{K}} \int_{\R
^d}
\1_{y+Z(s)}(x)\1_A(x)\h^n(\D x)\Lambda\bigl(\D(y,s)
\bigr)
\\
&=&\int_{\R^d\times\mathbf{K}}\int_{\R^d}\1_A(
\xi)\1_{Z(s)}(u)\lambda (\xi-u,s)\h^n(\D u)Q(\D s)\,\D\xi
\\
&=&\int_A \biggl(\underbrace{\int_{\mathbf{K}}
\int_{\R^d} \1_{Z(s)}(\xi-v)\lambda(v,s)
\h^n(\D v)Q(\D s)}_{=:\lambda_{\Theta_n}(\xi)} \biggr)\,\D\xi 
\end{eqnarray*}
and so the assertion.
\end{pf}

In \cite{AKV}, Proposition~9, it has been proved that for a class of
germ-grain models in $\R^d$ with independent and identically
distributed grains with finite $\h^n$-measure, $n<d$, the
probability that a point $x$ belongs to the intersection of two or
more enlarged grains is infinitesimally faster than $r^{d-n}$. The
i.i.d. assumption on the grains seems to be too restrictive; we
now extend it to more general germ-grain models as in above
assumptions. To this end, we shall make use of the assumption
(A3), which provides an integrability condition on the second
factorial moment measure $\nu_{[2]}$ of $\Phi$, similar to the
condition given on the intensity measure $\Lambda$ in (A2). Such a
result will be fundamental in the proof of the main theorem about
the validity of equation \eqref{aim}.
%
\begin{proposition}\label{prop}
Under the assumptions in Section \ref{secass},
the probability that a point
$x\in\R^d$ belongs to the intersection of two or more enlarged
grains $(y+Z(s))_{\oplus r}$ is infinitesimally faster than
$r^{d-n}$.

\end{proposition}

\begin{pf}
Let us observe that\vspace*{-1pt}
\begin{eqnarray*}
&&\E\biggl[\mathop{\sum_{%
(y_i,s_i),(y_j,s_j)\in\Phi,}}_{
y_i\neq y_j
} \1_{(y_i+Z(s_i))_{\oplus r}\cap
(y_j+Z(s_j))_{\oplus r}}(x)\biggr]
\\[-2pt]
&&\quad\stackrel{\eqref{nu2def}} {=}\int_{(\R^d\times\mathbf{K})^2}\1_{(x-Z(s_1))_{\oplus
r}}(y_1)
\1_{(x-Z(s_2))_{\oplus r}}(y_2)\nu_{[2]}(\D y_1,\D
s_1,\D y_2,\D s_2)
\\[-2pt]
&&\quad=\int_{\R^d\times\mathbf{K}^2} \biggl(\1_{(x-Z(s_2))_{\oplus
r}}(y_2)\int
_{(x-Z(s_1))_{\oplus r}}g(y_1,s_1,y_2,s_2)
\,\D y_1 \biggr) \,\D y_2Q_{[2]}(\D s,\D y).\vspace*{-2pt}
\end{eqnarray*}
By Theorem \ref{Minkgen} with $\mu=g(\cdot,s,y,t)\h^d$, together
with (A1) and (A3), it follows\vspace*{-1pt}
\begin{eqnarray*}
&&\lim_{r\da
0}\frac{1}{b_{d-n}r^{d-n}}\int_{(x-Z(s_1))_{\oplus
r}}g(y_1,s_1,y_2,s_2)
\,\D y_1
\\[-2pt]
&&\quad=\int_{x-Z(s_1)}g(y_1,s_1,y_2,s_2)
\h^n(\D y_1)\qquad \forall s_1,s_2\in
\mathbf{K}, \forall y_2\in\R^d,\vspace*{-1pt}
\end{eqnarray*}
and the
limit is finite being $g(\cdot,s_1,y_2,s_2)$ locally bounded by
(A3), and $\h^n(Z(s))<\infty$ for any $s\in\mathbf{K}$ by (A1). As
$Z(s)$ is lower dimensional for any $s\in\mathbf{K}$, it is clear
that\vspace*{-1pt}
\[
\lim_{r\da0}\1_{(x-Z(s_2))_{\oplus r}}(y_2)=0 \qquad\mbox{for }
\h^d\mbox{-a.e. } y_2\in\R^d\ \forall
s_2\in\mathbf{K},\vspace*{-2pt}
\]
thus\vspace*{-2pt}
\[
\lim_{r\da0}\frac{1}{b_{d-n}r^{d-n}}\1_{(x-Z(s_2))_{\oplus
r}}(y_2)\int
_{(x-Z(s_1))_{\oplus
r}}g(y_1,s_1,y_2,s_2)
\,\D y_1=0\vspace*{-1pt}
\]
{\spaceskip=0.19em plus 0.05em minus 0.02em for
$\h^d\mbox{-a.e. } y_2\in\R^d,  \forall s_1,s_2\in\mathbf{K}.$
Furthermore, by \eqref{condeta1}, (A1) and (A3) it follows that for
any $r\leq1$}\vspace*{-1pt}
\begin{eqnarray*}
&&\1_{(x-Z(s_2))_{\oplus
r}}(y_2)\frac{1}{b_{d-n}r^{d-n}}\int_{(x-Z(s_1))_{\oplus
r}}g(y_1,s_1,y_2,s_2)
\,\D y_1
\\[-2pt]
&&\quad\leq\1_{(x-Z(s_2))_{\oplus
1}}(y_2)\frac{\h^d(\Xi(s_1)_{\oplus
r})}{b_{d-n}r^{d-n}}\sup_{y_1\in(x-Z(s_1))_{\oplus
r}}g(y_1,s_1,y_2,s_2)
\\[-2pt]
&&\quad\leq\frac{2^n4^d b_d}{\gamma
b_{d-n}}\h^n\bigl(\Xi(s_1)\bigr)
\xi_{x,B_1(x)}(s_1,y_2,s_2) .
\end{eqnarray*}
By assumption (A3), we have that
\[
\int_{\R^d\times\mathbf{K}^2} \frac{2^n4^d b_d}{\gamma b_{d-n}}\h^n(
\Xi(s_1) \xi_{x,B_1(x)}(s_1,y_2,s_2)
\,\D y_2 Q_{[2](\D s,\D t)}<\infty,
\]
so the
dominated convergence theorem implies
%
\begin{equation}
\lim_{r\da0}\frac{\E[\sum_{(y_i,s_i),(y_j,s_j)\in\Phi,
y_i\neq y_j} \1_{(y_i+Z(s_i))_{\oplus r}\cap
(y_j+Z(s_j))_{\oplus r}}(x)]}{b_{d-n}r^{d-n}}=0. \label{lim1}
\end{equation}

Let $W_r$ be the random variable counting the number of pairs of
different enlarged grains of $\Theta_n$ which cover the point $x$:
%
\begin{equation}
W_r:=\#\bigl\{(i,j), i< j \dvt x\in\bigl(y_i+Z(s_i)
\bigr)_{\oplus r} \cap\bigl(y_j+Z(s_j)
\bigr)_{\oplus r} \bigr\}; \label{Wr}
\end{equation}
then
\[
W_r\leq\mathop{\sum_{
(y_i,s_i),(y_j,s_j)\in\Phi,}}_{
y_i\neq y_j} \1_{(y_i+Z(s_i))_{\oplus r}\cap
(y_j+Z(s_j))_{\oplus r}}(x),
\]
and so
\[
0\leq\lim_{r\da
0}\frac{\pr(W_r>0)}{b_{d-n}r^{d-n}}\leq\lim_{r\da
0}
\frac{\sum_{k=1}^\infty k\pr(W_r=k)}{b_{d-n}r^{d-n}}=\lim_{r\da
0}\frac{\E[W_r]}{b_{d-n}r^{d-n}}\stackrel{\eqref{lim1}} {
\leq} 0,
\]
and so the assertion.
\end{pf}

We are ready now to state and prove the main result of the
section.
%
\begin{theorem}\label{thlambda-lim} Under the assumptions in Section \ref
{secass},
%
\begin{equation}
\lim_{r\da0}\frac{\mathbb{P}(x\in\Theta_{n_{\oplus r}})}{b_{d-n}r^{d-n}}=\lambda_{\Theta_n}(x), \qquad\h^d
\mbox{-a.e. } x\in\R^d. \label{P-lambda}
\end{equation}
\end{theorem}

\begin{pf}
Let $Y_r$ be the random variable counting the number of enlarged
grains which cover the point $x$:
\[
Y_r:= \sum_{(y_i,s_i)\in\Phi}\1_{(y_i+Z(s_i))_{\oplus r}}(x),
\]
and $W_r$ be the random variable defined in \eqref{Wr}. By the
proof of Proposition \ref{prop}, we know that
\[
\pr(W_r>0)=\mathrm{o}\bigl(r^{d-n}\bigr) \quad\mbox{and}\quad
\E[W_r]=\mathrm{o}\bigl(r^{d-n}\bigr) ;
\]
thus,
noticing now that
\[
W_r=\cases{ %
0, &\quad $\mbox{if }Y=0,1,$\vspace*{2pt}
\cr
1, &\quad
$\mbox{if }Y=2,$\vspace*{2pt}
\cr
\pmatrix{Y_r
\cr
2}, &\quad $\mbox{if } Y
\geq 3,$}
\]
we get
\[
\pr(Y_r=2)=\pr(W_r=1)\leq \pr(W_r>0)=\mathrm{o}
\bigl(r^{d-n}\bigr)
\]
and
\[
0\leq\E[Y_r;Y_r\geq 3]\leq\E[W_r;Y_r
\geq3]\leq\E[W_r]=\mathrm{o}\bigl(r^{d-n}\bigr) ,
\]
which imply
%
\begin{eqnarray}\label{P-Ecard}
\lim_{r\da
0}\frac{\pr(x\in\Theta_{n_{\oplus r}})}{b_{d-n}r^{d-n}}&=& \lim_{r\da0}\frac{\pr(Y_r>0)}{b_{d-n}r^{d-n}}=
\lim_{r\da
0}\frac{\pr(Y_r=1)+\mathrm{o}(r^{d-n})}{b_{d-n}r^{d-n}}=\lim_{r\da
0}\frac{\E[Y_r]}{b_{d-n}r^{d-n}}
\nonumber
\\[-8pt]
\\[-8pt]
\nonumber
&\stackrel{\eqref {Campformula}} {=}&\lim_{r\da0}\frac{1}{b_{d-n}r^{d-n}}\int
_{
\mathbf{K}}\int_{(x-Z(s))_{\oplus r}}\lambda(y,s)\,\D yQ(\D s).
\end{eqnarray}

By Theorem \ref{Minkgen} with $\mu(\D y)=\lambda(y,s)\,\D y$, it
follows that
\[
\lim_{r\da
0}\frac{1}{b_{d-n}r^{d-n}}\int_{(x-Z(s))_{\oplus r}}\lambda(y,s)\,\D
y= \int_{x-Z(s)}\lambda(y,s)\h^n(\D y),
\]
besides, by observing
that
\begin{eqnarray*}
&&\frac{1}{b_{d-n}r^{d-n}}\int_{(x-Z(s))_{\oplus
r}}\lambda(y,s)\,\D y
\\
&&\quad\leq\frac{\h^d((Z(s))_{\oplus
r})}{b_{d-n}r^{d-n}} \sup_{y\in(x-Z(s))_{\oplus r}}\lambda(y,s) \stackrel{
\eqref{condeta1},\eqref{locbound1}} {\leq} \frac{2^n4^d b_d}{\gamma
b_{d-n}}\h^n
\bigl(\Xi(s)\bigr)\widetilde{\xi}_{B_2(x)}(s)\qquad \forall r<2,
\end{eqnarray*}
assumption (A2) and the dominated convergence theorem
imply
\[
\lim_{r\da
0}\frac{1}{b_{d-n}r^{d-n}}\int_\mathbf{K}\int
_{(x-Z(s))_{\oplus
r}}\lambda(y,s)\,\D yQ(\D s)= \int_{
\mathbf{K}}
\int_{x-Z(s)}\lambda(y,s)\h^n(\D y)Q(\D s),
\]
and so, by
\eqref{formulalambda},
%
\begin{equation}
\lambda_{\Theta_n}(x)= \lim_{r\da
0}\frac{1}{b_{d-n}r^{d-n}}\int
_\mathbf{K}\int_{(x-Z(s))_{\oplus
r}}\lambda(y,s)\,\D yQ(\D s)\qquad
\mbox{for $\h^d$-a.e. }x\in\R^d.\label{lambda}
\end{equation}
Finally, the assertion
follows:
\[
\lim_{r \da0}\frac{\mathbb{P} (x\in\Theta_{n_{\oplus
r}})}{b_{d-n}r^{d-n}} \stackrel{\eqref{P-Ecard},\eqref{lambda}}
{=}\lambda_{\Theta_n}(x)\qquad \mbox{for $\h^d$-a.e. }x\in
\R^d.
\]
\upqed\end{pf}

\subsection{Corollaries and remarks}
We point out that equations \eqref{formulalambda} and
\eqref{P-lambda} have been proved in \cite{MeandensityBool},
Theorem~3.13,
for a general class of Boolean models $\Theta_n$ with
intensity measure $\Lambda$ of the type $\Lambda(\D(x,s))=f(x)\,\D
xQ(\D s)$, and so with position-independent grains and typical
grain $Z_0$, by using the explicit form of the capacity functional
of $\Theta_n$. Actually, Proposition \ref{proplambda} and
Theorem \ref{thlambda-lim} generalize to Boolean models with
position-dependent grains, as stated in the following corollary,
under the assumptions (A1) and (A2) only, in accordance with
the above mentioned result in \cite{MeandensityBool}.

\begin{corollary}[(Particular case: Boolean models)]\label{corBoolean} If
$\Theta_n$ is a Boolean model with intensity measure
$\Lambda(\D(x,s))=\lambda(x,s)\,\D xQ(\D s)$, then all the results
stated in the above section hold under assumptions \textup{(A1)} and \textup{(A2)}.
\end{corollary}
\begin{pf}
It is enough to note that assumption (A3) is implied by (A1) and
(A2). Indeed, by \eqref{poisson} $g(\cdot,
s,y,t)=\lambda(\cdot,s)\lambda(y,t)$, so that $g(\cdot, s,y,t)$ is
locally bounded and $\h^n(\operatorname{ disc}(g(\cdot, s,y,t)))=0$ by
(A2), whereas \eqref{locbound} holds with
$\xi_{a,K}:=\widetilde{\xi}_{K}(s)\1_{(a-Z(t))_{\oplus
1}}(y)\lambda(y,t)$, by observing that
\begin{eqnarray*}
&&\int_{\R^d\times\mathbf{K}^2}\h^n\bigl(\Xi(s)\bigr)
\xi_{a,K}(s,y,t)\,\D yQ_{[2]}(\D s,\D t)
\\
&&\quad=\int_{\R^d\times\mathbf{K}} \1_{(a-Z(t))_{\oplus1}}(y)\lambda(y,t)\,\D yQ(\D t) \int
_{\mathbf{K}}\h^n\bigl(\Xi(s)\bigr)\widetilde{
\xi}_K(s)Q(\D s),
\end{eqnarray*}
with $\int_{\mathbf{K}}\h^n(\Xi(s))\widetilde{\xi}_K(s)Q(\D
s)<\infty$ by (A2), and
\begin{eqnarray*}
\int_{\R^d\times\mathbf{K}} \1_{(a-Z(t))_{\oplus
1}}(y)\lambda(y,t)\,\D yQ(\D t)&\leq&
\int_{\mathbf{K}}\h^d\bigl(\bigl(a-Z(t)
\bigr)_{\oplus1}\bigr)\sup_{y\in
(a-Z(t))_{\oplus1}}\lambda(y,t) Q(\D t)
\\
&\stackrel{\eqref{condizioneeta}, \mathrm{(A1)}} {\leq}&\int_{\mathbf{K}}
\frac{\h^n(\Xi(t))}{\gamma}2^n4^d b_d
\sup_{y\in(a-Z(t))_{\oplus1}}\lambda(y,t) Q(\D t)
\\
&\stackrel{\mathrm{(A2)}} {\leq}& \int_{\mathbf{K}}\frac{\h^n(\Xi(t))}{\gamma}2^n4^d
b_d \widetilde{\xi}_{B_1(a)}(t) Q(\D t)\stackrel{\mathrm{(A2)}} {<}
\infty.
\end{eqnarray*}
\upqed\end{pf}

\begin{remark}[(Independent marking)]\label{indmark}  If the point process
$\Phi$ is an independent marking of $\widetilde{\Phi}$, then the
two-point mark distribution $M_{x,y}(\D s,\D t)$ in \eqref{nu2gen} is
independent of $x$ and $y$, so that $M_{x,y}(\D s, \D
t)=Q_{[2]}(\D s, \D t)=Q(\D s)Q(\D t)$; accordingly,
$g(x,s,y,t)=\widetilde{g}(x,y)$. As a consequence, assumption (A3)
simplifies by replacing $g(x,s,y,t)$ with $\widetilde{g}(x,y)$. We
also recall that $\widetilde{g}(x,y)$ can be written in terms of
the so-called \textit{pair-correlation function} $\rho(x,y)$ in this
way:
\[
\widetilde{g}(x,y)=\rho(x,y)\lambda(x)\lambda(y) .
\]
Moreover, if in particular $\lambda$ and $g$ are bounded, say by
$c_1$ and $c_2$ in $\R$, respectively, then the finiteness of the
integral in assumptions (A2) and (A3) is trivially satisfied by
(A1), by taking $\widetilde{\xi}_K(s)\equiv c_1$ in \eqref
{locbound1} and $\xi_{a,K}(s,y,t):= c_2\1_{(a-Z(t))_{\oplus1}}(y)$ in
\eqref{locbound}, and noticing that\looseness=-1
\begin{eqnarray*}
&&\int_{\R^d\times\mathbf{K}^2}\h^n\bigl(\Xi(s)\bigr)
\xi_{a,K}(s,y,t)\,\D yQ_{[2]}(\D s,\D t)
\\
&&\quad\leq c_2\int_{\mathbf{K}}\frac{\h^n(\Xi(t))}{\gamma}2^n4^d
b_d Q(\D t)\int_{\mathbb{K}}\h^n\bigl(\Xi(s)
\bigr)Q(\D s)\stackrel{\mathrm{(A1)}} {<}\infty.
\end{eqnarray*}\looseness=0
\end{remark}

\begin{example}\label{exphi-tilde}Simple examples of point
processes $\widetilde{\Phi}$ having
bounded intensity $\lambda$ and second moment density $\widetilde{g}$,
are, for
instance, the binomial process of $m$ points in a compact region
$W\subset\R^d$ with $\h^d(W)>0$, and the Mat\`ern cluster
process (e.g., see \cite{CIME}). We remind that for the binomial
process we have
$\lambda(x)=m/\h^d(W)$ and
$\widetilde{g}(x,y)=m(m-1)/(\h^d(W))^2$; whereas for a Mat\`ern
cluster process in $\R^2$ in which the parent process is a
uniform Poisson process with intensity $\alpha$, and each cluster
consists of $N\sim\operatorname{ Poisson}(m)$ points independently and
uniformly distributed in the ball $B_r(x)$, where $x$ is the
centre of the cluster, we have $\lambda=m\alpha$, and
$\widetilde{g}(x,y)=\alpha^2 m^2+\alpha m^2\h^2(B_r(x)\cap
B_r(y))/(\pi^2 r^4)\leq\alpha^2 m^2+\alpha m^2/(\pi r^2)$.
Other examples of processes with bounded intensity and second
moment density are considered for instance in~\cite{prokesova}.
These, together with Example \ref{exsegment}, which gives an insight
into the validity of assumption~(A1), provide simple examples where all the assumptions (A1)--(A3)
hold.
\end{example}

\begin{example}We mention that an important case of random sets
of dimension 1 is given by the so-called
\textit{fibre processes} (e.g., see
\cite{Benes-Rataj}); they can taken as models in different fields,
as Biology (e.g., fibre systems in soils
\cite{Benes-Rataj}, Section~3.2.3) and Medicine (e.g., modelling
vessels in certain angiogenesis processes
\cite{KV-IAS,KV-Morale}), and it is clear that assuming
stationarity or that the fibres are the grains of a Boolean model
might be too restrictive in applications. Now, we have results for
studying also the more general case in which the fibres are not
independent of each other; note that assumptions (A1) and (A2) are
generally satisfied in applications: (A1) is trivial, since fibres
are usually assumed to be rectifiable (see\vadjust{\goodbreak} also Example \ref{exsegment}),
while (A2) and (A3) hold whenever~$\lambda$ are $g$ are, for
instance, bounded and continuous, as observed in the remark above.

Moreover, Proposition \ref{proplambda} applies and an explicit
expression for $\lambda_{\Theta_1}$ can be obtained in
terms of the intensity measure of the process. In particular, in
order to provide an explicit example, let us notice that in
\cite{MeandensityBool}, Example 2, an explicit formula for
$\lambda_{\Theta_1}$ is given for an inhomogeneous segment Boolean
model in $\R^2$, whose segments have random length and
orientation; the same assumptions on the intensity measure also
apply now to more general segment processes, not necessarily
Boolean models (e.g., with $\widetilde{\Phi}$ as in Example~\ref
{exphi-tilde}).\vspace*{-3pt}
\end{example}

\begin{remark}[(``one-grain'' random set)]\label{remone-grain}It is
worth noting that, as a very particular case of point process
$\widetilde{\Phi}$, we may consider the case in which
$\widetilde{\Phi}=\{X\}$, that is it is given by only one random
point $X$ in $\R^d$. Obviously, in this case $g\equiv0$, and only
assumptions (A1) and (A2) have to be satisfied for the validity of
all the results stated above. Even if this case might seem
trivial, actually it can be taken as a model for several real
applications, and it is of great interest, because it emerges
that whenever a random closed set $\Theta_n$ can be described by
a random point $X\in\R^d$ (not necessarily belonging to
$\Theta_n$, e.g., its centre if $\Theta_{d-1}$ is the
surface of a ball centred in $X$ with random radius $R$) and its
random ``shape'' $Z:=\Theta_n-X$, then we may provide sufficient
conditions on $\Theta_n$ such that our main result
\eqref{P-lambda} holds. Note that in this case $\Lambda(\D(x,s))$
represents the probability that the point $X$ is in the
infinitesimal region $\D x$ with mark in $\D s$. For instance, if
the ``shape'' does not depend on the position and $X$ is uniformly
distributed in a bounded region $W\subset\R^d$, then
$\Lambda(\D(x,s))=\D xQ(\D s)/\h^d(W)$. Then, it emerges that the
key assumption on the random closed set $\Theta_n$ which implies
\eqref{P-lambda} is the geometric regularity assumption (A1) on
its grains. As a matter of fact, (A1) can be seen as the
stochastic version of the condition~\eqref{condizioneeta} which
ensures the existence of the $n$-dimensional Minkowski content of
each grain, whereas~(A2) and (A3) are just technical assumptions;
in particular (A3) allows us to prove the statement of
Proposition \ref{prop} (in the Boolean case, it is already
contained in (A1) and (A2)).\vspace*{-3pt}
\end{remark}

Under the above assumptions, it follows in particular
that $\Theta_n$ admits the so-called \textit{local mean
$n$-dimensional Minkowski content}, which has been introduced in
\cite{AKV}; namely $\Theta_n$ is said to admit local mean
$n$-dimensional Minkowski content if the following limit exists
finite for any $A\subset\R^d$ such that $\E[\h^n(\Theta_n\cap
\partial A)]=0$
\[
\lim_{r\da0}\frac{\E[\h^d(\Theta_{n_{\oplus r}}\cap
A)]}{b_{d-n}r^{d-n}}=\E\bigl[\h^n(
\Theta_n\cap A)\bigr] .\vspace*{-3pt}
\]

\begin{proposition}If $\Theta_n$ satisfies assumptions \textup{(A1)} and \textup{(A2)},
then it
admits local mean $n$-dimensional Minkowski content.\vspace*{-3pt}
\end{proposition}
\begin{pf*}{Sketch of the proof}  It is sufficient to prove that $\Theta_n$
satisfies the hypotheses of Theorem 4 in~\cite{AKV}.

We already observed in Remark \ref{remradon} that
$\E[\mu_{\Theta_n}]$ is finite on bounded sets.

By proceeding along the same lines as in the proof of Theorem 3.9
in \cite{MeandensityBool} (here, by defining
$\widetilde{\Theta}(\omega):=\bigcup_{(x_i,s_i)\in\Phi(\omega
)}x_i+\Xi(s_i)$,
where $\Xi(s_i)\supseteq Z(s_i)$ as in (A1), and\vadjust{\goodbreak}
$\eta(\cdot):=\h^n(\widetilde{\Theta}(\omega)\cap W_{\oplus
2}\cap
\cdot)/\h^n(\widetilde{\Theta}(\omega)\cap W_{\oplus2})$), it is
easy to see that the hypotheses of the above mentioned theorem
are fulfilled with $Y:=\h^n(\widetilde{\Theta}\cap W_{\oplus
2})/\gamma$.\vspace*{-2pt}
\end{pf*}

\begin{remark}By Theorem \ref{thlambda-lim} and the above
proposition, the following chain of equalities holds, for any
$A\subset\R^d$ such that $\h^d(\partial A)=0$ (which implies
$\E[\h^n(\Theta_n\cap
\partial A)]=0$, being $\E[\mu_{\Theta_n}]\ll\h^d$):
\begin{eqnarray*}
\int_A \lim_{r\downarrow0} \frac{\mathbb{P}(x\in
{\Theta_{n_{\oplus r}}})}{b_{d-n}r^{d-n}} \,\mathrm{ d} x&=&
\int_A \lambda_{\Theta_n}(x)\,\D x=\E\bigl[
\h^n(\Theta_n\cap A)\bigr]
\\
&=&\lim_{r\da0}\frac{\E[\h^d(\Theta_{n_{\oplus r}}\cap
A)]}{b_{d-n}r^{d-n}}=\lim_{r\downarrow0}\int
_A \frac{\mathbb
{P}(x\in
{\Theta_{n_{\oplus r}}})}{b_{d-n}r^{d-n}} \,\mathrm{ d} x
\end{eqnarray*}
as if we might exchange limit and integral, answering to the open
problem raised in \cite{AKV}, Remark~8.\vspace*{-2pt}
\end{remark}

As mentioned in \cite{AKV}, several problems in real applications
are related to the estimation of the mean density of lower
dimensional inhomogeneous random sets (see also \cite{KV-IAS} and
reference therein); in particular, as a computer graphics
representation of lower dimensional sets in $\R^2$ is anyway
provided in terms of pixels, which can offer only a 2-D box
approximation of points in~$\R^2$, it might be useful to have
statistical estimators of the mean density $\lambda_{\Theta_n}$
based on the volume measure $\h^d$ of the Minkowski enlargement of
$\Theta_n$. To this end, a consistent and asymptotically unbiased
estimator $\widehat{\lambda}_{\Theta_n}(x)$ of
$\lambda_{\Theta_n}(x)$ has been introduced in \cite
{MeandensityBool}, based on equation~\eqref{P-lambda}, for a class of Boolean
models with typical grain $Z_0$. Having now proved that
\eqref{P-lambda} holds for more general random closed sets, that
is not only in stationary settings or for Boolean models, but also
for non-stationary germ-grains models whose grains are not assumed
to be independent each other, the same simple proof of
Proposition 6.1 in \cite{MeandensityBool} still applies, so that
we may state the following result.\vspace*{-2pt}
%
\begin{corollary}\label{corstimatore} Let
$\Theta_n$ satisfy the assumptions, and $\{\Theta^i_n\}_{i\in\N}$
be a sequence of random closed sets i.i.d. as $\Theta_n$; then the
estimator $\widehat{\lambda}^N_{\Theta_n}(x)$ of
$\lambda_{\Theta_n}(x)$ so defined
\[
\widehat{\lambda}_{\Theta_n}^N(x) := \frac{\sum_{i=1}^N\1_{\Theta_n^i\cap B_{R_N}(x)\neq\varnothing}}{N
b_{d-n}R_N^{d-n}}\vspace*{-2pt}
\]
is asymptotically unbiased and weakly
consistent for $\h^d$-a.e. $x\in\R^d$, if $R_N$ is such that
\[
\lim_{N\to\infty}R_N= 0\quad \mbox{and} \quad\lim_{N\to\infty}NR_N^{d-n}=
\infty.\vspace*{-2pt}
\]
\end{corollary}
%
\begin{remark}
$\widehat{\lambda}_{\Theta_n}^N(x)$ can be written also in terms
of the so-called \textit{empirical capacity functional} of
$\Theta_n$, which we recall to be defined as \cite{feng}
$\widehat{T}^N_{\Theta_n}(K):=\frac{1}{N}\sum_{i=1}^N
\1_{\Theta_n^i\cap K\neq\varnothing}$ for any compact
$K\subset\R^d$:
\[
\widehat{\lambda}_{\Theta_n}^N(x) :=\frac{\widehat{T}^N_{\Theta
_n}(B_{R_N}(x)) }{
b_{d-n}R_N^{d-n}}.\vadjust{\goodbreak}
\]
For a more detailed discussion on
$\widehat{\lambda}_{\Theta_n}^N(x)$ and related open problems, we
refer to \cite{MeandensityBool}, Section~6.
\end{remark}

\section{Mean surface density and spherical contact distribution}\label{sec4}
Let us now consider a random closed set $\Theta$ in $\R^d$, with
$\h^d(\Theta)>0$. A problem of interest is then the existence
(and which is its value) of the limit
\[
\sigma_{\Theta}(x):=\lim_{r\da0}\frac{\pr(x\in\Theta_{\oplus
r}\setminus\Theta)}{r} .
\]
The quantity $\sigma_\Theta(x)$ is usually called the \textit{specific area of $\Theta$} at point $x$, and it has been
introduced in \cite{matheron}, page~50. The name \textit{specific area}
comes from the fact that, under suitable regularity assumptions on
the boundary of $\Theta$ (e.g., when $\Theta$ has Lipschitz
boundary, or it is union of convex sets, etc.), $\sigma_\Theta(x)$
might coincide with the mean density $\lambda_{\partial
\Theta}(x)$ of $\partial\Theta$, that is the density of the
measure $\E[\mu_{\partial\Theta}]$ on $\R^d$. Moreover, it is
clearly related to the existence of the right partial derivative
at $r=0$ of the so-called \emph{local spherical contact
distribution function $H_\Theta$ of $\Theta$}, the function from
$\R_+ \times\R^d$ to $[0,1]$ so defined
%
\begin{equation}
H_\Theta(r,x):=\pr(x\in\Theta_{\oplus r} | x\notin\Theta).
\label{contact}
\end{equation}
We refer to \cite{MeandensityBool} and
\cite{ECS10} (and reference therein) for a more detailed
discussion on $\sigma_\Theta$; we point out that only results for
Boolean models with position-independent grains has been given
there, whereas in \cite{Hug-Last} general germ-grains models are
considered assuming that the grains are convex, so that results
and techniques from convex and integral geometry can be applied.
In this last mentioned paper, some formulae for contact
distributions and
mean densities of inhomogeneous germ-grain models are to
be taken in \emph{weak form} (e.g.,
\cite{Hug-Last}, Theorem 4.1), unless to add further suitable integrability
assumptions (e.g., in \cite{Hug-Last}, Remark 4.4, the existence of a
dominating integrable function is assumed). Nevertheless, the
assumption of
convexity of the grains in \cite{Hug-Last} seems to be too
restrictive in possible real applications, and it hides the fact
that $\sigma_\Theta$ may be differ from the mean boundary density
$\lambda_{\partial\Theta}$ of $\Theta$, as discussed in \cite
{MeandensityBool}. Indeed, we remind that the value of $\sigma_{\Theta
}$ is
strictly related to the value of the so-called \textit{mean outer
Minkowski content of $\Theta$} (and so of its
grains), which depends on the $\h^{d-1}$-measure of the set of the
boundary points of $\Theta$ where the \textit{$d$-dimensional density} of
$\Theta$
is 0 or 1 or 1/2 (e.g., see \cite{outerMin,MeandensityBool} for
more details on this
subject). In order to extend some results provided in \cite
{MeandensityBool} to general random closed sets, we briefly recall
basics on the
outer Minkowski content notion.

\subsection{$d$-dimensional densities and outer Minkowski content} Let
$A\in\mathcal{B}_{R^d}$; the quantity
$\mathcal{SM}(A)$ defined as (see \cite{AColeV})
\[
\mathcal{SM}(A):=\lim_{r\downarrow
0}\frac{\mathcal{H}^d(A_{\oplus r}\setminus A)}{r},
\]
provided
that the limit exists finite, is called \emph{outer Minkowski
content of $A$}.
Note that if $A$ is lower dimensional, then
$\mathcal{SM}(A)=2\mathcal{ M}^{d-1}(A)$, whereas if
$A$ is a $d$-dimensional set, closure of
its interior, then $A_{\oplus r}\setminus A$ coincides with the
outer Minkowski enlargement of $\partial A$ at distance $r$.

In
\cite{outerMin} two general classes of subsets of $\R^d$ which
admit outer Minkowski content has been introduced; in particular
we remind the definition of the so-called \textit{class
$\mathcal{O}$} and a related result.

\begin{definition}[(The class $\mathcal O$)]\label{classO} Let $\mathcal
O$ be the class
of Borel sets $A$ of $\R^d$ with countably $\h^{d-1}$-rectifiable
and bounded topological boundary, such that
\[
\eta\bigl(B_r(x)\bigr)\geq \gamma r^{d-1} \qquad\forall x\in
\partial A,\ \forall r\in(0,1)
\]
holds for some $\gamma>0$ and some probability measure $\eta$ in
$\mathbb{R}^d$ absolutely continuous with respect to
$\mathcal{H}^{d-1}$.
\end{definition}
The \emph{$d$-dimensional density}
(briefly, \emph{density}) of $A$ at a point $x\in\R^d$ is defined
by \cite{AFP}
\[
\delta_d(A,x):=\lim_{r\downarrow0}\frac{\mathcal{H}^d(A\cap
B_r(x))}{\h^d(B_r(x))},
\]
provided that the limit exists. It is clear that $\delta_d(A,x)$ equals
$1$ for all $x$ in the
interior of $A$, and 0 for all $x$ into the interior of the
complement set of $A$, whereas different values can be attained at
its boundary points. It is well known (e.g., see \cite{AFP}, Theorem 3.61) that if $\h^{d-1}(\partial A)<\infty$, then $A$ has
density either 0 or 1 or 1/2 at $\h^{d-1}$-almost every point of
its boundary. For every $t\in[0,1]$ and every
$\mathcal{H}^d$-measurable set $A\subset\mathbb{R}^d$ let
\[
A^t:=\bigl\{x\in\mathbb{R}^d\dvt \delta_d(A,x)=t
\bigr\}.
\]
The set of points $
\partial^*A:=\mathbb{R}^d\setminus(A^0\cup A^1)
$ where the density of $A$ is neither $0$ nor $1$ is called
\emph{essential boundary} of $A$. It is proved (e.g., see
\cite{AFP}) that all the sets $A^t$ are Borel sets, and that
$\mathcal{H}^{d-1}(\partial^*A\cap
B)=\mathcal{H}^{d-1}(A^{1/2}\cap B)$ for all
$B\in\mathcal{B}_{\R^d}$. It follows that for any $A$ with
$\h^{d-1}(A)<\infty$, it holds
%
\begin{equation}
\h^{d-1}(A)=\h^{d-1}\bigl(A^{1/2}\bigr)+
\h^{d-1}\bigl(A^0\cap\partial A\bigr)+\h^{d-1}
\bigl(A^1\cap\partial A\bigr).\label{misurabordo}
\end{equation}
As
Theorem \ref{teominkcont} gives general sufficient conditions on
the existence of the Minkowski content of a lower dimensional set,
as the following theorem gives similar general sufficient
conditions for the existence of the outer Minkowski content.
%
\begin{theorem}[(\cite{outerMin})]\label{teoouter} The class $\mathcal O$
is stable
under finite unions and any $A\in\mathcal O$ admits outer Minkowski
content, given by
\[
\mathcal{SM}(A)=\h^{d-1}\bigl(A^{1/2}\bigr)+2\h^{d-1}
\bigl(\partial A\cap A^0\bigr)= \h^{d-1}\bigl(\partial^* A
\bigr)+2\h^{d-1}\bigl(\partial A\cap A^0\bigr).
\]
\end{theorem}
Note
that $\mathcal{SM}(A)=\h^{d-1}(A)$ if $\h^{d-1}(\partial
\Theta\cap(\Theta^0\cup\Theta^1))=0$. A local version of the
outer Minkowski content is given in \cite{outerMin}, Proposition~4.13.

We also remind
that Theorem \ref{Minkgen} is a generalization of
Theorem \ref{teominkcont}; similarly, the next theorem might be
seen as a generalization of Theorem \ref{teoouter}.
%
\begin{theorem}[(\cite{MeandensityBool})]\label{Minkgen2}
Let $\mu$ be a positive measure in $\R^d$ absolutely
continuous with respect to $\h^d$ with locally bounded density
$f$, and
let $A$ belong to $\mathcal{ O}$. If
$\h^{d-1}(\operatorname{ disc}f)=0$, then
%
\begin{equation}
\lim_{r\da0}\frac{\mu(A_{\oplus r}\setminus
A)}{r}=\int_{\partial^*A} f(x)
\h^{d-1}(\D x)+2\int_{\partial A\cap A^0}f(x)\h^{d-1}(\D
x).\label {minkgener}
\end{equation}
\end{theorem}
%
\subsection{Specific area and mean surface density}\label{secsigma}
Let us consider a random closed set $\Theta$ in $\R^d$ with
$\h^d(\Theta)>0$, such that it might be represented as an
``one-grain'' random set by giving its random shape $Z$ and its
random location $y$, that is by giving a marked point process
$\Phi=(y,s)$ with $\pr(\Phi(\R^d\times\mathbf{K})>1)=0$, so that
\[
\Theta=y+Z(s)
\]
as discussed in Remark \ref{remone-grain}. For sake of simplicity,
let $Z$ be compact (the case in which $Z$ is locally compact might
be handled by introducing a suitable compact window containing the
point $x$ considered). Of course $\partial\Theta=x+\partial Z$,
and so the regularity properties of $\partial\Theta$ coincide
with the regularity properties of $\partial Z$. Let $\Phi$ have
intensity measure $\Lambda(\D(x,s))=\lambda(x,s)\,\D xQ(\D s)$ such
that
\begin{longlist}[(A1$'$)]
\item[(A1$'$)] for any $(y,s)\in\R^d\times\mathbf{K}$, $y+\partial Z(s)$
is a
countably $\h^{d-1}$-rectifiable and compact subset of $\R^d$,
such that there exists a closed set $\Xi(s)\supseteq\partial
Z(s)$ such that $\int_\mathbf{K}\h^{d-1}(\Xi(s))Q(\D s)<\infty$
and
\[
\h^{d-1}\bigl(\Xi(s)\cap B_r(x)\bigr)\geq\gamma
r^{d-1} \qquad\forall x\in\partial Z(s),\ \forall r\in(0,1)
\]
for some $\gamma>0$ independent on $y$ and $s$;
\item[(A2$'$)] for any
$s\in\mathbf{K} $, $\h^{d-1}(\operatorname{ disc}(\lambda(\cdot,s)))=0$ and
$\lambda(\cdot,s)$ is locally bounded
such that
for any compact $K\subset\R^d$
\[
\sup_{x\in K_{\oplus\operatorname{
diam}(Z(s))}}\lambda(x,s)\leq\widetilde{\xi}_K(s)
\]
for some
$\widetilde{\xi}_K(s)$ with
$\int_{\mathbf{K}}\h^{d-1}(\Xi(s))\widetilde{\xi}_K(s)Q(\D s)<\infty$.
\end{longlist}

Note that assumption (A1$'$) guarantees that
$Z$, and so $\Theta$, belongs to the class $\mathcal{O}$; in
particular it is easy to see that $\Theta$ satisfies the
hypotheses of Lemma 3.10 in \cite{MeandensityBool}, which
implies that $\Theta$ admits \emph{local mean outer Minkowski
content}, that is:
%
\begin{equation}
\lim_{r\da
0}\frac{\E[\mathcal{H}^d((\Theta_{\oplus r}\setminus\Theta)\cap
A)]}{r} =\mathbb{E}\bigl[\h^{d-1}\bigl(
\partial^*\Theta\cap A\bigr)\bigr]+2\E\bigl[\h^{d-1}\bigl(
\Theta^0\cap\partial\Theta\cap A\bigr)\bigr]\label{minkmean}
\end{equation}
for any Borel set $A$
with $\mathbb{E}[\h^{d-1}(\partial
\Theta\cap\partial A)]=0$ (and so for any $A$ with $\h^{d}(
\partial A)=0$, being $\E[\mu_{\partial\Theta}]\ll\h^d$).

The assumption (A2$'$) allows us to apply Theorem \ref{Minkgen2}
to prove that
\[
\sigma_{\Theta}(x)= \lambda_{\partial^*\Theta}(x)+ 2\lambda_{\Theta^0\cap\partial
\Theta}(x)
,\qquad \h^d\mbox{-a.e. }x\in\R^d,
\]
having denoted
by $\lambda_{\partial^*\Theta}$ and
$\lambda_{\Theta^0\cap\partial\Theta}$ the density of the measure
$\E[\h^{d-1}(\partial^*\Theta\cap \cdot)]$ and
$\E[\h^{d-1}(\Theta^0\cap\partial\Theta\cap \cdot)]$,
respectively; namely, we prove the following theorem.

\begin{theorem}\label{thed} Let $\Theta=y+Z$ be a random closed set
as above,
satisfying assumption \textup{(A1$'$)} and~\textup{(A2$'$)}; then
%
\begin{equation}
\sigma_\Theta(x):=\lim_{r\da0}\frac{\pr(x\in
\Theta_{\oplus
r}\setminus\Theta)}{r}=
\lambda_{\partial^*\Theta}(x)+ 2\lambda_{\Theta^0\cap\partial\Theta}(x), \qquad\h^d
\mbox{-a.e. }x\in\R^d.\label{sigma2}
\end{equation}
In particular, if
%
\begin{equation}
\int_{\mathbf{K}}\h^{d-1}\bigl(\partial^*Z(s)\bigr) Q(\D s)=
\int_{\mathbf{K}}\h^{d-1}\bigl(\partial Z(s)\bigr) Q(\D s)
\label {Pzeta},
\end{equation}
then
\[
\sigma_\Theta(x)=\lambda_{\partial\Theta}(x)=\int_{\mathbf
{K}}
\int_{x-\partial
Z(s)}\lambda(y,s)\h^{d-1} (\D y)Q(\D s),\qquad
\h^d\mbox{-a.e. }x\in\R^d.
\]
\end{theorem}

\begin{pf}
By applying the same arguments used in the proof of
Proposition 3.8 in \cite{MeandensityBool}, it follows that
$\E[\h^{d-1}(\partial\Theta\cap \cdot)]$ is absolutely
continuous with respect to $\h^d$ (and so $\E[\h^{d-1}(\partial^*
\Theta\cap \cdot)]$ and $\E[\h^{d-1}(\Theta^0\cap\partial
\Theta\cap \cdot)]$ as well, being $\partial^*\Theta$ and
$\Theta^0\cap\partial\Theta$ disjoint subsets of
$\partial\Theta$); the equation~\eqref{minkmean} is equivalent to
write
%
\begin{equation}
\lim_{r\da0}\int_A \frac{\pr(x\in\Theta_{\oplus
r}\setminus\Theta)}{r}=\int
_A \bigl(\lambda_{\partial^*\Theta}(x)+ 2\lambda_{\Theta^0\cap\partial\Theta}(x)
\bigr)\,\D x.\label{scambio}
\end{equation}
We want to apply the dominated convergence
theorem in order to exchange
limit and integral in the equation above.

Let us first prove that there exist the limit of $\pr(x\in
\Theta_{\oplus r}\setminus\Theta)/r$ for $r\da0$:
\begin{eqnarray*}
\lim_{r\da0} \frac{\pr(x\in\Theta_{\oplus
r}\setminus\Theta)}{r} &=& \lim_{r\da
0}\frac{\pr(\Phi\{(y,s) \dvt x\in(y+Z(s))_{\oplus r}\setminus
(y+Z(s))\}>0)}{r}
\\
&=& \lim_{r\da
0}\frac{\Lambda(\{(y,s) \dvt y\in(x-Z(s))_{\oplus r}\setminus
(x-Z(s))\})}{r}
\\
&=& \lim_{r\da
0}\frac{1}{r}\int_{\mathbf{K}}\int
_{(x-Z(s))_{\oplus
r}\setminus(x-Z(s))}\lambda(y,s)\,\D y Q(\D s).
\end{eqnarray*}
By applying now Theorem \ref{Minkgen2}, we get
\begin{eqnarray*}
&&\lim_{r\da0}\frac{1}{r} \int_{(x-Z(s))_{\oplus
r}\setminus(x-Z(s))}\lambda(y,s)
\,\D y
\\
&&\quad\stackrel{\eqref{minkgener}} {=}\int_{x-\partial^* Z(s)}\lambda(y,s)
\h^{d-1}(\D y)Q(\D s)+ 2 \int_{\mathbf{K}}\int
_{(x-\partial Z(s))\cap(x-
Z^0(s))}\lambda(y,s) \h^{d-1}(\D y),
\end{eqnarray*}
besides we
observe that
\begin{eqnarray*}
\frac{1}{r} \int_{(x-Z(s))_{\oplus
r}\setminus(x-Z(s))}\lambda(y,s)\,\D y&\leq&
\frac{1}{r} \int_{(x-\partial Z(s))_{\oplus r}}\lambda(y,s)\,\D y\\
&\leq&
\frac{\h^{d-1}(Z(s))}{r}\sup_{y\in(x-\partial Z(s))_{\oplus r}} \lambda(y,s)
\\
&\stackrel{\eqref{condeta1},\mathrm{(A2')}} {\leq}&
\h^{d-1}\bigl(\Xi (s)\bigr)\frac{2^{3d-1}
b_d}{\gamma}\widetilde{\xi}_{B_2(x)}(s)\qquad
\forall r<2.
\end{eqnarray*}
Therefore, assumption (A2$'$) and the dominated convergence theorem
imply
%
\begin{eqnarray}\label{star3}
\lim_{r\da0} \frac{\pr(x\in
\Theta_{\oplus r}\setminus\Theta)}{r}&=&\lim_{r\da
0}\frac{1}{r}\int
_{\mathbf{K}}\int_{(x-Z(s))_{\oplus
r}\setminus(x-Z(s))}\lambda(y,s)\,\D y Q(\D s)
\nonumber\\
&=&\int_\mathbf{K} \biggl(\int_{x-\partial^* Z(s)}
\lambda(y,s) \h^{d-1}(\D y)Q(\D s)\\
&&\hspace*{18pt}{}+ 2 \int_{\mathbf{K}}\int
_{(x-\partial
Z(s))\cap(x- Z^0(s))}\lambda(y,s) \h^{d-1}(\D y) \biggr) Q(\D
s).\nonumber
\end{eqnarray}
Analogously, for any fixed bounded Borel set $A$ and for any
$r<2$,
\[
\frac{\pr(x\in\Theta_{\oplus r}\setminus\Theta)}{r}\leq\int_{\mathbf{K}} \frac{\h^{d-1}(\Xi(s))}{\gamma}2^{3d-1}b_d
\widetilde{\xi}_{K}(s) Q(\D s)\stackrel{\mathrm{(A2')}}
{=}c\in\R,
\]
where $K$ is a compact subset of $\R^d$ containing $A_{\oplus
2}$.

Thus we may change limit and integral in \eqref{scambio}, and we
get
%
\begin{equation}
\lim_{r\da0}\int_A \frac{\pr(x\in\Theta_{\oplus
r}\setminus\Theta)}{r}=\int
_A \sigma_\Theta(x)\,\D x= \int_A
\bigl(\lambda_{\partial^*\Theta}(x)+ 2\lambda_{\Theta^0\cap
\partial
\Theta}(x) \bigr)\,\D x
\label{star2}
\end{equation}
for any $A$ with $\h^{d}(
\partial A)=0$, and so equation \eqref{sigma2} holds.

Assumption \eqref{Pzeta} ensures that the $\h^{d-1}$-measure of
the boundary of $Z$ equals the $\h^{d-1}$-measure of its essential
boundary, and so $\E[\h^{d-1}(\partial^*\Theta\cap
\cdot)]=\E[\h^{d-1}(\partial\Theta\cap \cdot)]$; in particular it
follows that $\lambda_{\partial^*\Xi}(x)=\lambda_{\partial\Xi}(x)$
and $ \lambda_{\partial\Xi\cap\Xi^0}(x)=0$ for
$\h^d\mbox{-a.e. } x\in\R^d$, and that
\begin{eqnarray*}
&&\int_{x-\partial^* Z(s)}\lambda(y,s) \h^{d-1}(\D y) + 2 \int
_{(x-\partial Z(s))\cap(x-
Z^0(s))}\lambda(y,s) \h^{d-1}(\D y)\\
&&\quad = \int
_{x-\partial
Z(s)}\lambda(y,s)\h^{d-1} (\D y).
\end{eqnarray*}
Thus, by \eqref{star3} and
\eqref{star2} we get
\[
\sigma_\Theta(x)=\lambda_{\partial\Theta}(x)=\int_{\mathbf
{K}}
\int_{x-\partial
Z(s)}\lambda(y,s)\h^{d-1} (\D y)Q(\D s),\qquad
\h^d\mbox{-a.e. }x\in\R^d.
\]
\upqed\end{pf}

\begin{remark}The above theorem answers also to the open problem posed
by Matheron in \cite{matheron}, page~50, about the equality between
the specific area $\sigma_{\Theta}$ and the mean boundary density
$\lambda_{\Theta}$ for a general random set $\Theta$. Again, such
an equality strongly depends on the geometric regularities of
$\partial\Theta$; of course the cases in which
$\sigma_{\Theta}\neq\lambda_{\partial\Theta}$ are, in a certain
sense, ``pathological,'' because condition \eqref{Pzeta} is
usually fulfilled in applications.
\end{remark}

Of course the specific area $\sigma_{\Theta}$ may be evaluated for
germ-grain processes whose grains have integer dimension $n<d$
($n=0$ is trivial), but it is clear that $\sigma_\Theta(x)\equiv
0$ if $n<d-1$.

In the case $d-1$, that is $Z(s)=\partial Z(s)$
for any $s\in\mathbf{K}$, assumptions (A1) and (A2) given in the
previous section coincide with (A1$'$) and (A2$'$) above; by noticing
that $\partial Z(s)=Z^0(s)\cap
\partial Z(s)$, and that $\pr(x\in\Theta)=0$ a.s., the results
\eqref{P-lambda} and \eqref{formulalambda} proved in
Theorem \ref{thlambda-lim} and Proposition \ref{proplambda},
respectively, are in accordance with Theorem \ref{thed}:
\begin{eqnarray*}
\sigma_{\Theta}(x)&=&\lim_{r\da0}\frac{\pr(x\in\Theta_{\oplus
r})}{r}=2
\lambda_{\Theta}(x)\\
&=&2\int_{\mathbf{K}}\int_{x-Z(s)}
\lambda (y,s)\h^{d-1}(\D y)Q(\D s).
\end{eqnarray*}

We point out that it seems to be hard to find out explicit
expressions for $\sigma_\Theta$ when $\Theta$ is a general
germ-grain model (i.e., non-Boolean) with $\h^d(\Theta)>0$, in
terms of its grains as we did for $\lambda_{\Theta_n}$ in
Proposition \ref{proplambda} in the $n$-dimensional case. Indeed,
due to the fact that the interior of the grains is in general not
empty, we cannot follow the same lines of the proof of the
mentioned proposition, because $\E[\h^{d-1}(\partial\Theta\cap
\cdot)]\neq\E[\sum_{(y_i,s_i)\in\Phi}\h^{d-1}((y_i+\partial
Z(s_i))\cap\cdot)]$.

Instead, when $\Theta$ is a Boolean model, and so thanks to
the independence property of its grains and to the knowledge of
the associated capacity functional, it
is possible to prove an explicit expression for its specific area,
as proved in \cite{MeandensityBool}, Proposition~3.7, in the case of
position-independent grains. By similar arguments of the previous
sections, it is easy to extend it to the case of a general Boolean
model $\Theta$ whose grains satisfy the above assumption (A1$'$)
and (A2$'$), obtaining that
%
\begin{eqnarray}
\label{sigmabool}&&\sigma_{\Theta}(x)=\pr(x\notin \Theta ) \biggl[ \int
_{\mathbf{K}}\int_{x-\partial^* Z(s)}\lambda(y,s)
\h^{d-1}(\D y)Q(\D s)
\nonumber
\\[-8pt]
\\[-8pt]
\nonumber
&&\hspace*{69pt}\quad{}+ 2 \int_{\mathbf{K}}\int_{(x-\partial Z(s))\cap(x-
Z^0(s))}\lambda(y,s)
\h^{d-1}(\D y)Q(\D s) \biggr].
\end{eqnarray}
We
may notice that the above expression for $\sigma_\Theta$ applies
only to Boolean models, thanks to independence properties of the
underlying point process $\Phi$, and that it cannot be true for
different germ-grain models: it is sufficient to consider the case
when $\Theta$ is an ``one-grain'' random set as in
Theorem \ref{thed}, and observe that its specific area given in
\eqref{star3} differs from \eqref{sigmabool}, being
$\pr(x\notin\Theta)\neq1$, in general.

\subsection{The spherical contact distribution function}
We are now able to give a general expression for the derivative in
$r=0$ of the spherical contact distribution function $H_{\Theta}$,
defined in \eqref{contact}, under the same general assumptions on
the random set $\Theta$ given in the previous section.

By noticing that $\pr(x\notin\Theta)H_{\Theta}(r,x)=\pr(x\in
\Theta_{\oplus r}\setminus\Theta)
$ and $H_{\Theta}(0,x)\equiv0$, the following corollary of
Theorem \ref{thed} is easily proved.
%
\begin{corollary}Let $\Theta$ be a random
closed set as in Theorem \ref{thed}; then
\begin{eqnarray*}
\frac{\partial}{\partial r}H_{\Theta}(r,x)_{|r=0}&=&\frac{\sigma_{\Theta
}(x)}{\pr(x\notin\Theta)}\\
&=&
\frac{\lambda_{\partial^*\Theta}(x)+
2\lambda_{\Theta^0\cap\partial\Theta}(x)}{\pr(x\notin\Theta)}, \qquad\h^d\mbox{-a.e.
}x\in\R^d,
\end{eqnarray*}
where the above derivative has
to be intended the right derivative in 0.

If in particular \eqref{Pzeta} is satisfied, then
\[
\frac{\partial}{\partial
r}H_{\Theta}(r,x)_{|r=0}=\frac{\lambda_{\partial\Theta}(x)}{\pr
(x\notin
\Theta)},\qquad
\h^d\mbox{-a.e. }x\in\R^d.
\]
\end{corollary}

\begin{remark}[(Boolean model and ``one-grain'' random set)] By the
corollary above and by \eqref{sigmabool} and \eqref{star3}, we
get the following explicit formulas in the case $\Theta$ is a
Boolean model (reobtaining~\cite{MeandensityBool}, equation~(4.1), as
particular case), or $\Theta$ is an ``one-grain'' random set:
\begin{eqnarray*}
&&\frac{\partial}{\partial
r}H_{\Theta}(r,x)_{|r=0}
\\
&&\quad =\cases{
\mbox{Boolean model}
\vspace*{2pt}\cr
\displaystyle\int_{\mathbf{K}}\int_{x-\partial^* Z(s)}\lambda (y,s)
\h^{d-1}(\D y)Q(\D s)\vspace*{2pt}\cr
\quad{}\displaystyle + 2 \int_{\mathbf{K}}\int
_{(x-\partial Z(s))\cap(x-
Z^0(s))}\lambda(y,s) \h^{d-1}(\D y)Q(\D s),
\vspace*{2pt}\cr
\mbox{``one-grain'' random set}
\vspace*{2pt}\cr
\displaystyle\frac
{ \int_{\mathbf{K}} [\int_{x-\partial^*
Z(s)}\lambda(y,s) \h^{d-1}(\D y)+ 2 \int_{(x-\partial Z(s))\cap
(x- Z^0(s))}\lambda(y,s) \h^{d-1}(\D y) ]Q(\D s)}{
\int_{\R^d\times\mathbf{K}}(1-\1_{x-Z(s)}(y))\lambda(y,s)\,\D yQ(\D
s) }. }
\end{eqnarray*}
\end{remark}

 In \cite{MeandensityBool}, Theorem~4.1, has been proved
a result about the differentiability of $H_{\Theta}$ with respect
to $r$ for a quite general class of Boolean models with typical
grain having \emph{positive reach}. Such a result can be easily
extended for Boolean models with position dependent grains by
considering an intensity measure $\Lambda(\D(y,s))$ of the type
$\lambda(y,s)\,\D yQ(\D s)$, instead of the type $f(y)\,\D yQ(\D s)$,
and by modifying the assumption of the cited theorem accordingly.
Here we reformulate such a result also for ``one-grain'' random
sets. In order to do this, we briefly recall some basic
definitions from geometric measure theory.

For any closed subset $A$ of $\mathbb{R}^d$, let $\operatorname{
Unp}(A):=\{x\in\mathbb{R}^d \dvtx \exists! a\in A \mbox{ such that } \operatorname{ dist}(x,A)=| a-x|\}.$
The definition of $\operatorname{ Unp}(A)$ implies the existence of a
projection mapping $\xi_A \dvt \operatorname{ Unp}(A)\to A$ which assigns to
$x\in\operatorname{ Unp}(A)$ the unique point $\xi_A(x)\in A$ such that
$\operatorname{ dist}(x,A)=| x-\xi_A(x)|$; then for all $x\in\operatorname{ Unp}(A)$
with $\operatorname{ dist}(x,A)>0$ we may define $u_A(a):=(x-\xi_A(x))/\operatorname{
dist}(x,A)$. The set of all $x\in\R^d\setminus A$ for which
$\xi_A(x)$ is not defined it is called \emph{exoskeleton of $A$},
and it is denoted by $\operatorname{ exo}(A)$. The \emph{normal bundle of
$A$} is the measurable subset of $\partial A\times\mathbf{S}^{d-1}$
defined by $N(A):=\{(\xi_A(x),u_A(x)) \dvt x\notin A\cup\operatorname{
exo}(A)\}. $ For any $x\in\partial^+ A:=\{x\in\partial A \dvt (x,u)\in N(A)\mbox{ for some }u\in\mathbf{S}^{d-1}\}$, we define
\[
N(A,x):=\bigl\{u\in\mathbf{S}^{d-1} \dvt (x,u)\in N(A)\bigr\}
\]
and
\[
\partial^1 A:=\bigl\{x\in\partial^+ A \dvt \operatorname{ card}N(A,x)=1\bigr\}.
\]
Note that for any $x\in\partial^1 A$, the unique element of
$N(A,x)$ is the \emph{outer normal of $A$ at $x$}, denoted here by
$n_x$. The reach of a compact set $A$ is defined by (see
\cite{Fed})
\[
\operatorname{ reach}(A):=\inf_{a\in A}\sup\bigl\{r>0 \dvt B_r(a)\subset
\operatorname{Unp}(A)\bigr\};
\]
for any set $A\subset\R^d$ with positive reach, the \emph{curvature
measures} $\Phi_{i}(A; \cdot )$ on $\R^d$, for $i=1,\ldots,
d-1$, introduced in \cite{Fed}, are well defined.

Then, by following the same lines of Section 4 in \cite
{MeandensityBool}, it is not difficult to prove the following
proposition for
an ``one-grain'' random set.
%
\begin{proposition}Let $\Theta$ be a random closed
set as in Theorem \ref{thed}, with $\operatorname{ reach}(Z(s))>R$ for some
$R>0$ and such that $\h^0(N(Z(s),x))=1$ for
$\h^{d-1}\mbox{-a.e. }x\in\partial Z(s),$ for all
$s\in\mathbf{K}$. Moreover, we assume that
\[
\int_{\mathbf{K}}|\Phi_i|\bigl(Z(s)\bigr) Q(\D s)<\infty\qquad
\forall i=1,\ldots,d-1,
\]
where $|\Phi_i|(Z)(s)$ is the total variation of
the measure $\Phi_i( Z(s);\cdot )$, and that the intensity
$\lambda(\cdot, s)$ is bounded, Lipschitz with Lipschitz constant
$\operatorname{ Lip} f(\cdot,s)$ such that $\int_{\mathbf{K}}\operatorname{ Lip}
f(\cdot,s)Q(\D s)<\infty$, and the set where $\lambda(\cdot,s)$ is
not differentiable is $\h^{d-1}$-negligible. Then, for all
$x\in\R^d$,
\begin{eqnarray*}
&&\frac{\partial}{\partial r} H_\Theta(r,x)=\frac{ \int_{\mathbf{K}}\int_{x-\partial
Z(s)_{\oplus r}}\lambda(y,s)\h^{d-1}(\D y) Q(\D s)}{
\int_{\R^d\times\mathbf{K}}(1-\1_{x-Z(s)}(y))\lambda(y,s)\,\D yQ(\D
s)}\qquad \forall r\in[0,
R),\label{Hprima}
\\
&&\frac{\partial^2}{\partial r^2}H(r,x)_{|r=0}\\
&&\quad=\frac{ \int_{\mathbf{K}} [2\pi\int_{\R^d}
\lambda(y,s)\Phi_{d-2}(x-Z(s);\D
y) + \int_{x-\partial^1Z(s)}D_{n_y}\lambda(y,s)\h^{d-1}(\D
y) ]Q(\D s)}{
\int_{\R^d\times\mathbf{K}}(1-\1_{x-Z(s)}(y))\lambda(y,s)\,\D yQ(\D
s)},
\end{eqnarray*}
where $D_{n_y}\lambda(\cdot,s)$ is the
directional derivative
of $\lambda(\cdot, s)$ along $n_y\in\mathbf{S}^{d-1}$.
\end{proposition}



\section*{Acknowledgements} I wish to thank Prof. Luigi Ambrosio (SNS,
Pisa), Prof. Vincenzo Capasso (University of Milan) and the
anonymous referees for their careful reading of the paper and for
their useful comments and suggestions.

%

%


\printhistory

\end{document}